\documentclass[12pt,twoside, reqno]{amsart}
\usepackage[utf8]{inputenc}
 \usepackage[english]{babel}
\usepackage{amsthm}
\usepackage{amssymb, verbatim}
\usepackage{amsmath}
\usepackage{eucal}
\usepackage{mathrsfs}
\usepackage{hyperref}

 \topmargin = - 30 pt 
\usepackage[all]{xy}
\usepackage{amsmath,amsthm,amssymb,amsopn,amscd,tikz-cd}

\newtheorem{thm}{Theorem}[section]

\newtheorem{prop}[thm]{Proposition}

\theoremstyle{definition}

\newtheorem{defn}[thm]{Definition}

\newtheorem{ex}[thm]{Example}

\newtheorem{rem}[thm]{Remark}

\numberwithin{equation}{section}

\begin{document}

\title[Uniqueness of braided tensor structures]{Uniqueness of braided tensor structures and the Huang--Lepowsky structure on affine vertex operator algebras  in the classical Lie types and $G_2$}

 \author[C.~Pinzari]{Claudia Pinzari}

 \email{pinzari@mat.uniroma1.it}
\address{Dipartimento di Matematica, Sapienza Universit\`a  di
Roma\\ P.le Aldo Moro, 5 -- 00185 Rome, Italy}\maketitle

  \centerline{\it Dedicated to Sergio Doplicher, with gratitude}

\begin{abstract}

We prove that the unitary modular tensor category structure that we constructed
for all Lie types in our previous work, on the module category of the affine
vertex operator algebra $V_{\mathfrak{g}_{k}}$ at positive integer level,
coincides with the Huang--Lepowsky braided tensor structure for all the
classical Lie types and $G_2$. This answers, in the vertex operator algebra
setting and for these types, the problem of proving the Finkelberg equivalence
theorem directly, without the equivalence of Kazhdan--Lusztig. Simply-laced and
non-simply-laced types are treated uniformly throughout: our framework rests on
Wenzl's work on the unitary structure of the quantum group fusion categories,
which applies to all simple Lie algebras, whereas the original work of Kazhdan
and Lusztig explicitly considered the simply-laced case.
 
The proof uses  two abstract uniqueness results, established here with
complete proofs: in a semisimple pre-tensor category with a generating object
satisfying a duality property with respect to a braiding, the associativity and
braiding morphisms are determined by their restrictions to a small collection of
objects. The duality property is then verified for the classical types and $G_2$
by appealing to known generalized quantum Schur--Weyl duality.
 
This completes our previous work, where the global quantum gauge group
$A_{W}(\mathfrak{g},q)$, an intrinsic construction for braided tensor
$C^*$-categories in the sense of the Doplicher--Roberts duality problem, was
canonically constructed and used to endow the Zhu algebra of
$V_{\mathfrak{g}_{k}}$ with what we called a unitary coboundary weak quasi-Hopf
structure with $3$-coboundary associator, via an isometric analytic Drinfeld
twist along Wenzl's continuous de-quantization curve, for all Lie types. The
appeal to the braid group  duality, made here for the first time in this context, is
what unifies the Doplicher--Roberts construction of a compact gauge group in
higher-dimensional algebraic quantum field theory with the Huang--Lepowsky
braided tensor structure and  quantum gauge groups.
 
\end{abstract}

    \tableofcontents

\section{Introduction}

The main result of this paper is that Zhu's linear equivalence between the
module category of the affine vertex operator algebra $V_{\mathfrak{g}_k}$ at
positive integer level, carrying the Huang--Lepowsky modular tensor category
structure, and the representation category of its Zhu algebra, carrying the
structure induced by the unitary coboundary weak quasi-Hopf structure of the
latter defined and constructed in our previous work for all Lie types, is an equivalence of ribbon braided tensor categories for the Lie types
$A$, $B$, $C$, $D$ and $G_2$. Composed with the equivalences recalled below,
this identifies ${\rm Rep}_{\rm HL}(V_{\mathfrak{g}_k})$ with the quantum group
fusion category $\mathcal{C}(\mathfrak{g},q)$ for those types. In the vertex
operator algebra setting, and for these types, this proves the Finkelberg
equivalence theorem directly, without the equivalence of Kazhdan--Lusztig, and
so answers the problem raised by Y.-Z.~Huang \cite{Huang2018}. The precise
statement, with the quantum parameter and the twist made explicit, is
Theorem~\ref{new} in Section~\ref{Main}.
 
The argument treats simply-laced and non-simply-laced Lie algebras uniformly.
This is not incidental: our framework rests on Wenzl's work on the unitary structure of the quantum group fusion category, which applies to
all simple Lie algebras, whereas the original work of Kazhdan and Lusztig
explicitly considered the simply-laced case, and explicit statements extending
their construction to the non-simply-laced case are not easy to isolate in the
literature. In types $B$, $C$ and $G_2$ our results therefore cover positive
integer levels that the original argument does not reach.
 
The technical core of the paper, and what carries the result above, is a pair of
abstract uniqueness theorems, Theorem~\ref{claim0} and its generalization Theorem~\ref{claim1}.
They assert that in a semisimple pre-tensor category with a generating object
satisfying a duality property with respect to a braiding, the associativity and
braiding morphisms are determined by their restrictions to a small collection of
objects. The duality property is then verified for the classical types and $G_2$
by a generalized quantum Schur--Weyl duality, in Section~\ref{12}. This uses a
property from the discovery of quantum groups that we had not exploited in
\cite{CGP}: the extension of Schur--Weyl duality between the braid group and the
quantum group inside the centralizer algebras of the fusion powers of the
fundamental representation. The statements and proofs of
Theorem~\ref{claim0} and Theorem~\ref{claim1}, and the verification of
Section~\ref{12}, are self-contained: they depend on no result quoted from
elsewhere, and may be read and checked independently of the rest of the paper.
 
Three results are quoted from our earlier work \cite{CGP} and used in the proof
of Theorem~\ref{new}: the construction of the unitary coboundary weak Hopf
algebra $A_W(\mathfrak{g},q)$, whose representation category is equivalent to
$\mathcal{C}(\mathfrak{g},q)$ as a unitary ribbon braided tensor category
(Theorem~2.5 there, and Theorem~9.2 of \cite{CP} in type $A$); the twist
isomorphism between the Zhu algebra and a twist of $A_W(\mathfrak{g},q)$
(Theorem~2.8 there); and the resulting equivalence with the quantum group fusion
category (Corollary~2.9 there). These are stated in full in
Section~\ref{4}, so that the logical status of every ingredient of
Theorem~\ref{new} is visible without consulting \cite{CGP}. Nothing else is
assumed from that work, the remaining material being established here:
the intervening sections  set out the construction in the form the present argument requires, and include results that do not appear there.

 

For the exceptional types $E$ and $F$ the duality property with respect to the
braid group representation is not known at present; this is precisely what the
generalized quantum Schur--Weyl duality of Section~\ref{12} supplies in the
other types. The results of our earlier work nonetheless apply: First,
${\rm Rep}(V_{\mathfrak{g}_k})$ carries a unitary rigid ribbon braided tensor
structure, lifted from the unitary coboundary structure of the Zhu algebra and
equivalent to $\mathcal{C}(\mathfrak{g},q)$, so that the
Finkelberg--Kazhdan--Lusztig equivalence holds unconditionally with respect to
this transported structure; see Remark~\ref{Exceptional_types}. Second, the
comparison with the Huang--Lepowsky structure already goes a considerable way in
these types: by Theorem~\ref{Finkelberg_HL}\,(b) and (c), the two structures
have the same tensor product bifunctor, Zhu's equivalence satisfies the ribbon
equations for all objects, and the braiding and associativity equations hold for
the pairs
$(V_\lambda, V)$, $(V, V_\lambda)$
and for the triples
$(V_\lambda, V, V)$, $(V, V_\lambda, V)$, $(V, V, V_\lambda)$,
with $V_\lambda$ an arbitrary irreducible object of the fusion category. What
remains open is the propagation of these identities to all objects.
 
Both results belong to a wider program, which addresses two problems at the
intersection of algebraic quantum field theory, the representation theory of
vertex operator algebras, and quantum groups:
\begin{itemize}
\item the problem posed by Doplicher and Roberts on the intrinsic construction
of global quantum gauge groups for braided tensor $C^*$-categories \cite{DR};
\item the problem posed by Huang of proving the Finkelberg equivalence theorem
directly, avoiding the recourse to the Kazhdan--Lusztig equivalence at negative
shifted levels \cite{Huang2018}.
\end{itemize}
In \cite{CGP} we constructed, for the WZW model and all Lie types, a unitary
modular tensor category structure on the module category of the affine vertex
operator algebra at positive integer level, together with a global quantum gauge
group for that structure; the present paper identifies that structure with the
Huang--Lepowsky one in the classical types and $G_2$, and thereby completes the
picture in those cases.
 
Our approach obtains rigidity directly from the quantum group fusion category.
It bypasses the negative-level shift of the original argument, the Verlinde
formula and the monodromy of the Knizhnik--Zamolodchikov equations ordinarily
required in the vertex operator algebra setting, and the Jones index used for
conformal nets \cite{Wassermann}; rigidity is thus established prior to global
modular properties, entirely within the vertex operator algebra framework.
 
In the rest of the introduction we explain these results in further detail. 

\subsection{Historical context and  motivations}
The relationship between the representation categories of affine Lie algebras and quantum groups has historically been mediated by the monodromy of the Knizhnik--Zamolodchikov (KZ) equations. This   approach, however, encounters a primary difficulty at positive integer levels due to singularities in the KZ equations.

Major work was carried out in a series of papers by Kazhdan and Lusztig
   in the mid 90s  \cite{KLseries}.   
The authors  constructed the structure of a rigid braided tensor category  on  a   category of
modules  ${\mathcal O}_\ell$ of affine Lie algebras associated to simply laced simple Lie algebras ${\mathfrak g}$ at complex levels $k$ such that $\ell=k+\check{h}\notin{\mathbb Q}_{\geq0}$.  Moreover,  they proved that ${\mathcal O}_\ell$ is braided tensor equivalent to a non-semisimple category of modules of the Drinfeld-Jimbo quantum group
 $U_q({\mathfrak g})$ at the  root of unity $q=e^{i\pi/\ell}$.

As detailed in Appendix~\ref{rigidity_problem}, at positive integer levels the rigorous construction of   braided  tensor categories was given by Huang and Lepowsky in a series of works in a very general   setting of vertex operator algebras, and rigidity and modularity were first proved by  Y.-Z. Huang. The proof of rigidity
relies on the Verlinde formula.

 Finkelberg  constructed a braided tensor equivalence from a category $\tilde{\mathcal O}_\ell$ of modules of affine Lie algebras at positive integer levels $k$  introduced by Beilinson, Feigin and Mazur \cite{BFM} to
a subquotient category  $\tilde{\mathcal O}_{-\ell}$
of Kazhdan-Lusztig category using duality between affine Lie algebra modules   \cite{Finkelberg}, \cite{Finkelberg_erratum}. Invoking   the equivalence by Kazhdan and Lusztig, he established a braided tensor equivalence $\tilde{\mathcal O}_\ell\to{\mathcal C}({\mathfrak g}, q)$ to the corresponding subquotient of the quantum group category. His proof of rigidity of the category relies on the use of the Verlinde formula as well.

The Kazhdan--Lusztig rigid tensor category construction does not apply to certain negative rational shifted levels. Consequently,  Finkelberg's proof of rigidity of  $\tilde{\mathcal O}_{\ell}$ can not  cover some exceptional cases--namely, $E_6$ and $E_7$ at level $k=1$, and $E_8$ at levels $k=1,2$.  
Furthermore,
 explicit statements extending the Kazhdan--Lusztig rigid braided tensor category construction to non-simply laced Lie algebras are not easy to isolate in the literature, owing to  complexities of the corresponding KZ equations.
 In addition, as explicitly noted by Finkelberg (see  Sect. 2.6 in \cite{Finkelberg} and part 2 (f) in \cite{Finkelberg_erratum}), the rigidity of $\tilde{\mathcal O}_{\ell}$ and  the tensor equivalence property of the Finkelberg functor may   be unproven for   other non-simply laced Lie algebras at low levels ($k\leq2$).

Conversely, in the operator algebraic approach to conformal field theory (CFT), rigid braided tensor $C^*$-categories are constructed directly, and the Verlinde formula is derived as a consequence of categorical traces. A similar procedure applies to quantum group fusion categories; for the Drinfeld category, Drinfeld derived rigidity by constructing an explicit isomorphism and a twist between $U_h(\mathfrak{g})$ and his quasi-Hopf algebra \cite{Drinfeld_quasi_hopf}. Because a twist operation preserves rigidity, he   established the rigidity of the Drinfeld category. 
These distinct approaches form the core motivation for our work (see Section~\ref{2new} and Appendix~\ref{2}  for further details).

Our framework overcomes the aforementioned limitations    by constructing a direct equivalence with the quantum group fusion category that holds uniformly for both simply laced and non-simply laced Lie algebras at all positive integer levels. By doing so, this approach completely bypasses the   recourse
 to the KZ equations, the negative-level shifting, and the traditional reliance on the Verlinde formula to prove rigidity. Furthermore, our proof of rigidity avoids the use of the Jones index within the von Neumann algebra framework of the conformal net WZW model  \cite{Wassermann}. Consequently, this framework restores the natural categorical framework---establishing   rigidity prior to global modular properties---entirely within the vertex operator algebra setting. 

Our main method is the construction of an independent {\it unitary} rigid ribbon braided tensor category in ${\rm Rep}(V_{{\mathfrak g}_k})$ 
derived from ${\mathcal C}({\mathfrak g}, q)$
via the introduction of a {\it unitary structure}, and construction of semisimple quantum gauge groups inspired by algebraic quantum field theory. Our equivalence with the quantum group fusion category is established in Corollary 2.9 in \cite{CGP}. This result builds on the work by Wenzl \cite{Wenzl} and is inspired by the work of Drinfeld \cite{Drinfeld_quasi_hopf}.

The aim of this paper is to compare our ribbon braided tensor category structure on ${\rm Rep}(V_{{\mathfrak g}_k})$ established in \cite{CGP} with the Huang-Lepowsky ribbon braided tensor category structure ${\rm Rep}_{\rm HL}(V_{{\mathfrak g}_k})$. This comparison is stated in Theorem \ref{new} and is presented in a more complete form, combined with our previous results, in Theorem \ref{Finkelberg_HL}. We emphasize that our quantum group and tensor category approach yields a   uniform treatment of both simply laced and non-simply laced Lie algebras--covering types $A, B, C, D,$ and $G_2$ for the Huang-Lepowsky structure, and all Lie types for our analytic structure--valid   for all   positive integer levels.

\subsection{A unified framework: connecting abstract duality and vertex operator algebras}
We approach these problems by unifying the Doplicher--Roberts abstract duality program with the VOA framework. An abstract duality theory is a form of Tannaka--Krein duality where, rather than being equipped with a predefined fibre functor from the outset, one must construct the fibre functor intrinsically from the category itself. By resuming the study of Mack--Schomerus quantum symmetries \cite{MS} , we explicitly unify the treatment of high-dimensional algebraic quantum field theory (AQFT) with the low-dimensional, braided setting. We provide a canonical construction of Mack--Schomerus symmetry algebras that is   compatible with the Huang--Lepowsky ribbon braided tensor structure for affine vertex operator algebras at positive integer levels.

The conceptual bridge between the Doplicher--Roberts program and the VOA setting is established by utilizing the {\it initial terms of primary fields} as a unifying link. Since a primary field is uniquely determined by its initial term  \cite{Tsuchiya_Kanie}, we exploit this fact to connect the fusion tensor product of $\mathcal{C}(\mathfrak{g},q)$ to the fusion of the corresponding representations of the vertex operator algebra. 

In concretely setting up this connection, we build on two  milestones:
\begin{itemize}
    \item On the quantum group side, we employ the work of Wenzl on the unitary structure of $\mathcal{C}(\mathfrak{g},q)$ \cite{Wenzl}, based on the use of a fundamental representation of $U_q(\mathfrak{g})$ and the construction of the unitary fusion tensor addendum submodule with any other irreducible object in the fusion category.
    \item On the conformal field theory side, the connection via initial terms of primary fields is enlightened by the work of Wassermann  \cite{Wassermann} on the loop group conformal net of ${\rm SU}(N)$.
\end{itemize}
This allows us to realize abstract duality within the concrete setting of the affine VOA under the Huang--Lepowsky tensor structure.

To develop the canonical construction of Mack--Schomerus symmetry algebras systematically, we connect the quantum group and VOA settings via two natural linear faithful functors. On the quantum group side, we employ the Wenzl functor on the unitary fusion category $\mathcal{C}(\mathfrak{g},q)$. On the VOA side, we employ the linear minimum energy functor (Zhu's functor) on the affine module category $\text{Rep}(V_{\mathfrak{g}_k})$, which yields a linear equivalence with the module category of the Zhu algebra under standard semisimplicity conditions  \cite{Zhu}. The fusion of a fundamental representation of $U_q(\mathfrak{g})$ with any other irreducible in the fusion category analyzed in \cite{Wenzl} is in this way grounded into the fusion of the corresponding fundamental representations of the Zhu algebra with its respective irreducible counterparts.

\subsection{Core method: weak Hopf algebras and analytic Drinfeld twists}
In our previous work \cite{CGP}, we introduced the notion of a weak Hopf algebra as the Hopf counterpart to the Mack--Schomerus weak quasi-Hopf algebra structure. The algebras $A_W(\mathfrak{g},q)$---introduced in \cite{CP} for the type $A$ case and in \cite{CGP} for the remaining Lie types---are examples of weak Hopf algebras associated with the Wenzl functor. These structures, termed {\it unitary coboundary weak Hopf $C^*$-algebras}, possess a $C^*$-structure and a unitary structure on tensor products of representations defined by a positive operator induced by the Drinfeld coboundary matrix $\overline{R}$, making the underlying braid group representations spacially unitary.

By utilizing the unitary coboundary weak Hopf algebra $A_W(\mathfrak{g},q)$ and the associative Zhu algebra $A(V_{\mathfrak{g}_k})$, we transport the  rigid, modular, and unitary tensor structures from the quantum group fusion category directly to the affine vertex operator algebra. This is achieved by constructing an analytic Drinfeld twist $T$ for $A_W(\mathfrak{g},q)$ and an isomorphism to the Zhu algebra induced by Wenzl's continuous quantization curve \cite{Wenzl}. Under this framework of weak quasi-Hopf algebras, the associator of the Zhu algebra is necessarily realized as a 3-coboundary. This enables us to define finite-dimensional tensor structures on the Zhu algebra from the twist and continuous curve of $A_W(\mathfrak{g},q)$, and then lift them to infinite-dimensional tensor structures on the vertex operator algebra.

We prove in Section~\ref{8} that Zhu's linear equivalence and its inverse transport all the categorical structure on the infinite-dimensional modules of the VOA endowed with the Huang--Lepowsky modular tensor structure \cite{Huang(modularity), Huang2} to the Zhu algebra. Thus, the structure on the VOA can be recovered completely from the transported representation category of the finite-dimensional Zhu algebra $A(V_{\mathfrak{g}_k})$. Building on the fusion of the fundamental representation with any other irreducible and Wenzl's   projections, we ultimately show that these projections extend to a fusion tensor product in a coherent way for the module category of the Zhu algebra, compatibly with the Huang--Lepowsky structure.

This approach yields a self-contained unitary rigid ribbon braided tensor structure on the module category $\text{Rep}(V_{\mathfrak{g}_k})$ of the affine VOA for all Lie types directly from the structure of the corresponding quantum group fusion category, along with a ribbon braided tensor equivalence. In particular, we obtain unitarity for all Lie types in a uniform way on $\text{Rep}(V_{\mathfrak{g}_k})$ because the same holds for the quantum group fusion category \cite{Wenzl}. This work was strongly influenced by a remark by Vaughan Jones on the property of manifest unitarity of unitary structures in CFT, which originally led us to connect VOAs from quantum groups by seeking such a unitary structure on the Zhu algebra. Consequently, the fusion rules, modular $S$-matrix, and the Verlinde formula emerge naturally as derived consequences rather than prior assumptions.

While our prior work \cite{CGP} established this independent unitary rigid ribbon braided tensor structure on $\text{Rep}(V_{\mathfrak{g}_k})$ lifted from the Zhu algebra and its equivalence with the modular fusion category of the quantum group, the present paper completes the program. We provide a complete identification of our analytic structure with the Huang--Lepowsky braided tensor structure for the classical and $G_2$ Lie types. To achieve this, we exploit a key property in the discovery of quantum groups that we had not previously considered in \cite{CGP}: the extension of Schur--Weyl duality between the braid group and the quantum group within the centralizer algebras of the fusion powers of the fundamental representation. Together, these two papers allow us to establish our main result.

\subsection{Geometric and analytic intuition of the twist}
Before stating the main theorem, we outline the geometric and analytic intuition underlying our twist framework. Our approach provides an alternative realization of the Huang--Lepowsky rigid ribbon braided tensor structure on $\text{Rep}(V_{\mathfrak{g}_k})$ via the explicit construction of a unitary structure on this category. 

Specifically, the twist $T$ is a non-commutative analogue of a Radon--Nikodym derivative constructed from the square root of the self-adjoint coboundary of $U_q(\mathfrak{g})$ and the fusion tensor product of $\mathcal{C}(\mathfrak{g},q)$. While the positive coboundary matrix $\overline{R}$ of $A_W(\mathfrak{g},q)$ functions as a density operator defining the deformation of the inner product---carrying the information of the change in 'measure' from the standard inner product necessary to make the braiding unitary---the twist $T$ provides the explicit implementation of this change, functioning essentially as an isometry between the two structures. 

Together with the algebra isomorphism between $A_W(\mathfrak{g},q)$ and $A(V_{\mathfrak{g}_k})$ induced by Wenzl's de-quantization curve \cite{Wenzl}---which interpolates continuously from certain modules of the quantum group $U_q(\mathfrak{g})$ to modules of the classical compact group $G$---the twist isomorphism successfully ``flattens'' the quantum braided structure of the $A_W(\mathfrak{g},q)$-modules into the standard braided structure of the Zhu algebra representations following the two unitary structures. This illuminates conceptually why the VOA and quantum group categories are unitarily equivalent. From an algebraic point of view, our approach describes the Huang--Lepowsky associativity  morphisms of $\text{Rep}_{\rm HL}(V_{\mathfrak{g}_k})$ directly as a 3-coboundary of the twist $T$ on the Zhu algebra.

\subsection{Main Theorem}\label{Main}
Let $\mathfrak{g}$ be a complex simple Lie algebra and $k$ a positive integer. Let $\ell = k + h^{\vee}$, where $h^{\vee}$ is the dual Coxeter number, and let $d\in\{1,2,3\}$ be the ratio of the square lengths of the long to the short root. We define the quantum parameter $q=e^{\frac{i\pi}{d\ell}}$.

Let us consider the unitary coboundary weak Hopf algebra $A_{W}(\mathfrak{g},q)$ established in Theorem~2.5 of \cite{CGP} (Theorem~9.2 of \cite{CP} in the type $A$ case), which satisfies
\begin{equation}
    \text{Rep}(A_{W}(\mathfrak{g},q))\simeq\mathcal{C}(\mathfrak{g},q)
\end{equation}
as unitary ribbon braided tensor categories. Let the Zhu algebra $A(V_{\mathfrak{g}_k})$ be endowed with the unitary coboundary weak quasi-Hopf algebra structure derived from the explicit twist isomorphism:
\begin{equation}\label{Zhu}
    A(V_{\mathfrak{g}_k})\simeq(A_{W}(\mathfrak{g},q))_{T}
\end{equation}
established in Theorem~2.8 of \cite{CGP}. Let us endow $\text{Rep}(A(V_{\mathfrak{g}_k}))$ with the unitary modular fusion category (UMFC) structure induced by this unitary coboundary weak quasi-Hopf algebra on $A(V_{\mathfrak{g}_k})$. We also consider the module category $\text{Rep}_{\rm HL}(V_{\mathfrak{g}_k})$ of the affine vertex operator algebra $V_{\mathfrak{g}_k}$ endowed with the   Huang--Lepowsky modular tensor category structure \cite{Huang(modularity), Huang_2, Huang2}.

The following is the main result of the present paper.

\begin{thm}\label{new}
Zhu's linear equivalence
\begin{equation}
    {\rm Rep}_{\rm HL}(V_{\mathfrak{g}_k})\rightarrow {\rm Rep}(A(V_{\mathfrak{g}_k}))
\end{equation}
becomes a ribbon braided tensor equivalence for the Lie types A, B, C, D, and $G_{2}$, realizing $\text{Rep}_{\rm HL}(V_{\mathfrak{g}_k})$ as a UMFC. Hence, we obtain the complete chain of ribbon braided tensor equivalences:
\begin{equation}
    {\rm Rep}_{\rm HL}(V_{\mathfrak{g}_k})\simeq {\rm Rep}(A(V_{\mathfrak{g}_k}))\simeq {\rm Rep}(A_{W}(\mathfrak{g},q))\simeq\mathcal{C}(\mathfrak{g},q).
\end{equation}
\end{thm}

\begin{rem}[The exceptional types and  perspective]
\label{rem:exceptional-outlook}
The   realization of the exceptional Lie types $E$ and $F$ within the explicit generator-lineage framework relies upon the current state of generalized quantum Schur--Weyl duality. Beyond the specific settings where such explicit duality with the braid group holds, a   perspective exists for treating   examples   where the braid group is not generating. Namely, the intrinsic reconstruction of global quantum symmetries might be handled--independently of explicit centralizer generation properties--via the general machinery of the braided Doplicher--Roberts program \cite{pinzari_survey}.
\end{rem}

\begin{rem}[Organization of the Proof]
\label{rem:proof-navigation}
The technical and intermediate steps required to complete the proof of Theorem~\ref{new} are    formulated in Theorem~\ref{Finkelberg_HL} below.
\end{rem}

\subsection{Organization of the paper}
The paper is organized as follows.
 
\text{Section~\ref{2new}} states the problems posed by Huang and by
Doplicher--Roberts, and sets out our global quantum gauge group strategy for the
Finkelberg--Kazhdan--Lusztig theorem.
 
\text{Section~\ref{4}} gives the definition of a unitary coboundary weak Hopf
algebra and recalls, in the form used here, the construction of
$A_W(\mathfrak{g},q)$ from \cite{CGP}.
 
\text{Section~\ref{5}} states the main result, Theorem~\ref{new}, together with
the technical Theorem~\ref{Finkelberg_HL}, parts (a) to (d).
 
\text{Section~\ref{6}} sets out the strategy of the proof and identifies, step
by step, which ingredients are quoted from \cite{CGP} and which are established
in the present paper.
 
\text{Section~\ref{7}} sets out the elements of Tannaka--Krein reconstruction
associated with a fibre functor possessing a weak quasi-tensor structure.
 
\text{Section~\ref{8}} proves part (a) of Theorem~\ref{Finkelberg_HL}. We
describe the tensor structure of the module category of the Zhu algebra obtained
by transporting the Huang--Lepowsky tensor structure on elementary building
blocks, and we embed these tensor products into the full ordinary tensor product
space as a canonical addendum submodule compatible with the transport
construction. Proposition~\ref{transport} and the remarks following it are not
in \cite{CGP}; they are stated here in the form the argument requires.

\text{Section~\ref{9}} establishes the principal unitarity properties of the
categorical transport. It recalls the Drinfeld twist and the Wenzl continuous
curve implementations realizing the twist isomorphism between
$A_W(\mathfrak{g},q)$ and the Zhu algebra, which makes the latter a unitary
coboundary weak quasi-Hopf algebra with $3$-coboundary associator, and shows
that the transported braiding and associativity morphisms coincide with those
obtained from the twist method. The last part of the section introduces the
cohomological reduction tool, via the pentagon and hexagon equations, that
extends the equivalence from the building blocks to all tensor products of
objects; this is carried out in Section~\ref{11}.

\text{Section~\ref{10}} sets out the structural analogies and divergences
between the original Drinfeld--Kohno theorem and our analogue from \cite{CGP},
and shows how the unitary structure permits the bypass of the KZ equations in
constructing the ribbon braided tensor structure on the module category of the
vertex operator algebra.

\text{Section~\ref{11}} contains the main new results of the paper. We
state and prove Theorem~\ref{claim0} and its generalization
Theorem~\ref{claim1}, which reduce the coincidence of two braided tensor
structures to a duality property with respect to the braiding. Their proofs are
self-contained and independent of \cite{CGP}.

\text{Section~\ref{12}} completes the proof of part (d) of
Theorem~\ref{Finkelberg_HL}, by verifying the duality property and centralizer
algebra generation, through generalized quantum Schur--Weyl duality, for the
classical Lie types and $G_2$.
 
\text{Section~\ref{13}} is devoted to concluding remarks and outlook.
 
\text{Appendix~\ref{2}} provides historical background on the interplay between
primary fields and quantum groups.
 
\text{Appendix~\ref{rigidity_problem}} reviews the historical development of
rigid braided tensor categories for affine Lie algebras and vertex operator
algebras, the proofs of rigidity by Huang--Lepowsky in the vertex operator
algebra setting and by Finkelberg in the affine Lie algebra setting, and the
underlying rigidity problem.

 \section{Doplicher-Roberts and Huang problems} \label{2new}

 While the equivalence between quantum group fusion categories and categories of modules of affine Lie algebras at positive integer level   was proposed in the mid 1990s, the complete
 proof of rigidity of  $\tilde{\mathcal O}_{\ell}$ and derivation of the tensor equivalence was given in 2013  using the Verlinde formula
  \cite{Finkelberg_erratum} building on work by
 Faltings and Teleman.
   We refer the reader to the appendix for more historical remarks concerning the proof 
  of rigidity for affine vertex operator algebras and affine Lie algebras at positive integer levels.
  
   The ribbon braided tensor  equivalence established in \cite{Finkelberg}, \cite{Finkelberg_erratum}   relies on a passage through negative levels and then applies the equivalence with the quantum group fusion category by Kazhdan and Lusztig. 
 As a consequence, Finkelberg's proof excludes certain exceptional cases like $\mathfrak{g}=E_8$ at level $k=2$ since this case is excluded in   \cite{KLseries}. 
 This led us to consider the following problem posed by  Y.-Z. Huang   in \cite{Huang2018}.

\medskip

\noindent{\bf Problem 4.4 in  \cite{Huang2018}.} {\it Find a direct construction of this   equivalence without using the equivalence given by Kazhdan-Lusztig so that this equivalence covers all the cases, including the important ${\mathfrak g}=E_8$ and $k=2$ case.}\medskip

Moreover, in light of Drinfeld's construction of the ribbon braided tensor equivalence in Drinfeld-Kohno theorem,
and his consequent derivation of rigidity of the Drinfeld category, is natural to  ask whether the Verlinde formula is necessary to prove rigidity, and whether rigidity is necessary to prove the ribbon braided tensor equivalence $\tilde{\mathcal O}_\ell \to {\mathcal C}({\mathfrak g}, q)$.

We also mention   Gannon's  observation in Sect. 6.2.3 of \cite{Gannon_book} regarding the
desire of finding a direct connection  between quantum groups and CFT.
Our work provides a direct constructive method--addressing Huang's Problem 4.4 in \cite{Huang2018}--by utilizing a global quantum gauge group method originating in algebraic quantum field theory, and the Zhu algebra to establish the equivalence.

From the side of mathematical axiomatizations of CFT, 
  Huang and Lepowsky have done extensive work on the construction of tensor structures of modules and showing 
  all the properties of modular tensor category \cite{Huang_LepowskiI}, \cite{Huang_LepowskiII}, \cite{Huang_LepowskiIII}, \cite{Huang1}, \cite{Huang_Lepowski_affine}, \cite{Huang(modularity)}, \cite{Huang3}, \cite{Huang4},  \cite{Huang_2}, \cite{Huang2}.

An explicit construction of the equivalence is helpful
 to   see how  rigidity,
modularity or unitarity in the setting of CFT may   be transported from corresponding results in the setting of quantum groups at roots of unity, which seem simpler,   by our explicit equivalence.

\subsection{Inspiration from Algebraic Quantum Field Theory}
 Our starting point  to approach this problem    arose from a related problem   posed by Doplicher and Roberts     in the axiomatic setting of the operator algebraic approach to quantum field theory (AQFT).  \medskip
 
 \noindent{\bf Problem posed in Section 7 in   \cite{DR}.}
 {\it Can we extend the duality theory for compact groups [61] Theorem 6.1 to more general objects, by replacing strict symmetric monoidal C*-categories by strict braided monoidal C*-categories? What is the full class of compact group-like objects which arise in this way?}
 \medskip
 
 We next explain  this problem and connect it to Huang's problem.  
 Algebraic quantum field theory  was initiated by Haag and Kastler in  \cite{Haag_Kastler} in the 1960s, and developed
by   founding breakthrough works
by Doplicher, Haag, and Roberts \cite{DHR}, \cite{DHR2}, \cite{DR_qft1},  \cite{DR_qft2} in the 1970s. The operator algebraic approach to QFT 
joins high rigour provided by the operator algebra methods with first principles in physics on which the theory is formulated. We refer
to the   book by Haag \cite{Haag},  for an in-depth introduction to this theory and to the work by Mueger \cite{Mueger_DR}, \cite{Mueger_modular}, \cite{Mueger_duality}.

  In the context of high dimensional quantum field theories, categorical structures
  first emerged from the discovery   of the notion of localized endomorphism of the local observable algebra, allowing
composition of two of them which was understood as a strict tensor product structure.
The mentioned works by 
 Doplicher, Haag, and Roberts led to  a deep   construction of the structure of  a rigid symmetric strict tensor $C^*$-category with simple unit.

 Rigidity, called charge conjugation in this setting, describes particle-antiparticle duality. 
In their work rigidity is derived as a weak limit of inner automorphisms induced by transfer operators moving
the support of the given localized endomorphism $\rho$ to (space-like) infinity. This gives a left inverse $\Phi_\rho$ of $\rho$,
which is a completely positive map. From this the conjugate localized endomorphism $\overline{\rho}$ is derived.
The left inverse may be described in the category as a map induced by the adjoint action of a   solution of the conjugate equations on a   intertwiner translated to the right by the tensor product.
It  was    shown that the values of the statistics parameter
(which is the constant $\Phi_\rho(\varepsilon)$, with $\varepsilon$ the action of the generator of the permutation group
in the algebra of intertwining operators of the tensor powers of $\rho$), can be completely
classified, and it can only take the values given by the reciprocal of a positive
  integer $d$ (the order of the parastatistics) up to a sign   which classifies    charges of para-Bose
or para-Fermi type.
This analysis also needed to prove  that the statistics parameter can not take the value $0$ and was 
excluded by physically motivated conditions \cite{DR_qft2}, \cite{Haag}.

 Because of the representation of the permutation group in the category, together with    their results on the classification of the statistics parameter,    and the $C^*$-structure, in the 70s  it arose the question of whether that was the category of representations of a compact group, because those data arising from a single localised endomorphism were also satisfied by the representation category of a compact Lie subgroup of ${\rm SU}(d)$,  that would be a   gauge group of global symmetries. 
They  considered the abelian case   in the 70s \cite{DR_qft1}.
  In that   case, every object had dimension $1$, and in modern terminology these categories are referred to as pointed in e.g. \cite{EGNO}.

 In a series of     papers dated the end of the 80s,
  Doplicher and Roberts were able to  solve   the  problem in full generality.
Their construction gives 
  a unique compact gauge group with representation category equivalent to that
 of   a given symmetric rigid strict tensor $C^*$-category with simple unit, yielding the formulation of a new     abstract    duality theory for compact groups \cite{DR_duals}, \cite{DR1},\cite{DR_annals}, and led to   desired applications to high dimensional Algebraic QFT   \cite{DR}.

 The term {\it abstract duality} refers to the property that the given categories do not have a tensor functor to the category of f.d. vector spaces from the beginning, raising the problem of constructing these structures. It is thus
an abstract form of duality in the sense of Tannaka-Krein, where in the latter the tensor structure on the embedding functor is given from the start.
In their theory, an important role is played by the notion of symmetric functor to the category of f.d. vector spaces, for which they show existence and uniqueness.
This property allows to find their unique compact gauge group $G$.

\subsection{Quantum symmetries and primary field reconstruction in low-dimensional
AQFT} In low-dimensional quantum field theories (specifically in two and three spacetime dimensions), the categories of superselection sectors naturally admit a braided tensor structure rather than a symmetric one. Consequently, the standard  DHR Bose-Fermi alternative does not hold any longer, allowing for the existence of non-integer statistics. Early models exhibiting such behavior were first constructed by Streater and Wilde. The DHR theory of superselection sectors was subsequently generalized to accommodate braid group statistics by Fredenhagen, Rehren, and Schroer \cite{FRS, FRS2}, who showed that the representations of the braid group within these tensor $C^*$-categories are   unitary.

Motivated by the Drinfeld-Kohno theorem and the representation theory of quantum groups at roots of unity--where Hermitian representations of $U_q(\mathfrak{sl}_2)$ become non-semisimple and exhibit null vectors--Mack and Schomerus sought to recover the truncated fusion tensor products of conformal field theory (CFT) from a semisimple quotient algebra. To this end, they extended Drinfeld's framework to define \emph{weak quasi-Hopf algebras}. They proposed that the true global quantum symmetries of minimal conformal models are weak quasi-Hopf algebras canonically associated with $U_q(\mathfrak{sl}_2)$ at roots of unity. As a first main example, they rigorously showed that the primary field operators of the chiral Ising model satisfy weak quasi-Hopf-covariant operator product expansions and obey local braid relations \cite{MS_quasi_quantum_group_symmetry}. Their program culminated in Schomerus's generalized construction of a weak quasi-Hopf algebra and a corresponding field algebra for any rational conformal field theory,   extending the Doplicher-Roberts reconstruction theorem to the braided case \cite{MS}, \cite{Schomerus}.

Mack-Schomerus  quantum symmetry program saw  further but limited development. A primary obstacle was that the emerging quantum gauge groups in the general theory were highly non-unique and possessed a complicated algebraic structure, largely due to the non-triviality of the associator. These ambiguities remained unclarified, leaving an incomplete reconciliation with the standard Doplicher-Roberts theory. Furthermore, the pioneering examples provided by Mack and Schomerus were   derived   from quantum group representation theory, rather than reconstructed intrinsically from the local observables. This approach raises  doubts as to whether their symmetries can be understood as weak versions of Hopf algebras--a notion and a connection   they did not address.

The algebraic results of Tsuchiya and Kanie were   developed into the operator algebraic approach to conformal field theory
 framework by 
A. Wassermann \cite{Wassermann}.
Utilizing von Neumann algebras, Wassermann constructed the algebras of observables
for chiral model of the CFT describing the WZW model, known as the   {\it loop group
conformal net}.

Subsequently, the theory developed toward the direct computation and classification   of the tensor categories of localized endomorphisms. By largely leaving the   quantum   gauge group viewpoint aside, this structural shift proved   successful
for 2D classifications. However, it also resulted in a   lack of a unifying mathematical language 
between low-dimensional braided theories and higher-dimensional symmetric theories governed by standard compact gauge groups.

.
\medskip

\subsection{Our Approach to the two Problems: Quantum Gauge Groups and Drinfeld Twists} 
 Our approach to Huang's problem and to Doplicher-Roberts problems is derived from synthesizing the 
   developments in affine Lie algebras by Wassermann and Toledano-Laredo on the loop group conformal net
    \cite{Toledano_laredo}
  with the unitary properties of quantum group fusion categories previously highlighted by Wenzl \cite{Wenzl}, from the view point 
  of quantum symmetries in algebraic conformal field theory.
  In our approach primary fields play a prominent role to establish the equivalence, as originally envisaged by Mack and Schomerus.
  In our approach, a weak quasi-Hopf algebra is naturally associated to the conformal field theory, while  a weak Hopf algebra,
  in a suitable new sense, is naturally associated to the quantum group fusion category, and we connect them by a Drinfeld twist.
  
Specifically, we   solve Huang's problem and advance Mack-Schomerus quantum symmetry program, by introducing a global {\it quantum gauge group} approach. We construct finite-dimensional weak Hopf $C^*$-algebras $A_W({\mathfrak g}, q)$ with unitary coboundary structures associated to the fusion category of quantum groups at roots of unity. To achieve this, we utilize a specific linear fiber functor $W$ originating from the work of Wenzl on
 ${\mathcal C}({\mathfrak g}, q)$ to the category of Hilbert spaces
\cite{Wenzl}.
 By applying an explicit Drinfeld twist--analogous to Drinfeld's original proof of the Drinfeld-Kohno theorem--we implement a unitary weak quasi-tensor structure directly onto the Zhu algebra. This explicit construction allows us to naturally and directly transport the rigid, modular, and unitary tensor structures from the quantum group category to the representation category of the affine vertex operator algebra (VOA) at positive integer levels, thereby providing the desired direct link in Huang's problem.

Our work has roots in our previous works in the setting of compact quantum groups of Woronowicz, \cite{CDPR}, \cite{DPiR}, \cite{P_braided}--\cite{PRinduction};
 however the different nature of categories involved here necessitates a very different framework.
Woronowicz developed   compact quantum groups as noncommutative function algebras in the setting of $C^*$-algebras 
via a concrete Tannaka-Krein duality for $C^*$-algebras extending the classical case. In his case, categories are naturally endowed with a fibre functor to the category of f.d. Hilbert spaces.

Our approach extends  Woronowicz Tannaka-Krein duality on several aspects. 
The first aspect is that our categories do not have a fibre functor a priori. Similarly to   Doplicher-Roberts setting, we must  construct one  naturally,  compatible with the braided tensor structure of the category. We need to do this first for the quantum group fusion category ${\mathcal C}({\mathfrak g}, q)$, and then for   the VOA fusion category ${\rm Rep}(V_{{\mathfrak g}_k})$.

To describe  our approach,  we go back to the constructive approach of rigid braided tensor categories motivated by CFT. We recall that Drinfeld original motivation
for introducing quasi-Hopf algebras was to interpret them as   symmetries   to describe braided tensor categories  from important models of rational conformal 
field theories motivated by the work by Moore and Seiberg. Although as noted in \cite{Drinfeld_exposition}, they could not describe precisely simple models such as the Ising, Drinfeld evaluated quasi-Hopf algebras  a good description of quantum symmetries for  the WZW model in an approximate way, i.e.
in the limit to infinity of the central charge.

Our approach to Huang problem    is  based on  the program of
 understanding Finkelberg-Kazhdan-Lusztig combined equivalence   as
an extension to quantum groups of Doplicher-Roberts   duality theory for the   class
of WZW models that are described by   affine Lie algebras at positive integer levels.  

From this perspective, we regard the equivalence described by  Drinfeld-Kohno theorem  as a 
Tannakian counterpart that needs to be implemented on the categories
appearing on the sides
of Kazhdan-Lusztig-Finkelberg equivalence
$\tilde{\mathcal O}_\ell\to{\mathcal C}({\mathfrak g}, q)$.     
In building this abstract duality theory, we were inspired by Neshveyev and 
Tuset's simplified proof of the Drinfeld-Kohno theorem \cite{NT_KL}  in the setting of compact quantum groups, 
viewing it as a Tannakian counterpart  implemented at the level of Hopf and quasi-Hopf algebras. 
The construction of the weak Hopf algebras, weak quasi-Hopf algebras starting from the category and a linear fibre functors
are not present in their work. Drinfeld original approach was the opposite, to start from the construction of Hopf and quasi-Hopf algebras.
Their methods for constructing the twist are inspired by Kazhdan-Lusztig work,   our
 methods  are inspired by Drinfeld original twist.
 Differences between the two situations arise from the different unitary structures
  at roots of unity or at real values of $q$.

The weak Hopf algebras $A_W({\mathfrak g}, q)$ that we   construct in \cite{CGP} from the quantum group fusion category are  generalizations of duals of the compact quantum groups of Wornowicz: They have  algebraic and unitarity properties that extend properly that of the   duals of the compact quantum groups. 

The main algebraic difference is that the coproduct $\Delta$ can not be unital in our case due to the fusion rules of the quantum grpoup fusion category. It is  also  not strictly associative either, but this lack of strict associativity is rather mild, because we construct a fibre functor satisfying two tensorial diagrams. Moreover the antipode satisfies the usual Hopf algebra axioms and commutes with the $^*$-involution.

The $^*$-involution of $A_W({\mathfrak g}, q)$ is anticomultiplicative, rather than 
multiplicative as in the case of compact quantum groups.  
 The hermitian structure of quantum group fusion categories   at roots of unity was   analysed in \cite{Wenzl}.
 This structure implies in a non-trivial way that the adjoint   bialgebra and the opposite bialgebra
 of  $A_W({\mathfrak g}, q)$ coincide, reflecting a similar structure of the hermitian structure
 of $U_q({\mathfrak g}, q)$--although $A_W({\mathfrak g}, q)$
 is not obtained by a co-ideal.
 To form tensor product of $^*$-representations, we need a 
 twist relating the adjoint   bialgebra of  $A_W({\mathfrak g}, q)$ with the given bialgebra. The existence of the twist allows to unify and extend the theory with that of the compact quantum groups. Our twist is inspired by the   twist that Drinfeld used for   a different purpose: the proof of the Drinfeld-Kohno theorem.
 To construct this twist, Drinfeld   made use of Drinfeld {\it coblundary $\overline{R}$-matrix}.
 In the case of symmetric tensor categories, the coboundary matrix is the identity. Otherwise it induces  an intertwiner
 $\rho\otimes\sigma\to\sigma\otimes\tau$ in the fusion category, similarly to how an $R$-matrix induces the braiding as an intertwiner.
 However for the coboundary, the intertwiner does not satisfy the braid group relation, but has   always of period $2$, extending the property of the permutation symmetry in the classical case.
We summarize the structure of $A_W({\mathfrak g}, q)$ by saying that this algebra  is a {\it unitary coboundary weak Hopf algebra}.

This is where the present paper closes a gap left open above. Our earlier 
construction produced an equivalence with the quantum group fusion category, but
 left it unrelated to the Huang--Lepowsky structure, which arises from a 
 different circle of ideas. In the Doplicher--Roberts theory it is the classical Schur--Weyl duality,
  between the permutation group and a compact group inside the intertwiner algebras of the tensor powers of a single object, that yields the uniqueness of 
  the symmetric fibre functor and hence the gauge group. 
  
  The duality we use here is its braided analogue: between the braid group and the quantum group, inside the centralizer algebras of the fusion powers of the fundamental representation. What the braided setting requires in addition has no counterpart in the symmetric case. There the permutation symmetric braiding is determined by its value on a single pair, the relations of the permutation group leaving no freedom, and no uniqueness statement is called for. Here the weaker relations of the braid group make it a theorem that, assuming a duality property with respect to a braiding, the associativity and braiding morphisms are determined by their restrictions to a small collection of objects; this is the content of Theorems~\ref{claim0} and \ref{claim1}.   It is the two together, the concrete generalized quantum Schur--Weyl duality and the abstract uniqueness results, that supply the unifying language whose absence we noted above.

 \section{Unitary coboundary weak Hopf algebras $A_W({\mathfrak g}, q)$}\label{4}
 The notion of {\it weak quasi-Hopf algebra} was first introduced by Mack and Schomerus motivated by the work of Drinfeld on quasi-Hopf algebras, specifically the Drinfeld-Kohno theorem  \cite{MS}.

 \subsection{Weak quasi-Hopf algebras}
 The weakness   refers to a non-unital coproduct, and is needed to construct a coalgebra structure
 on the quotient $A$ of $U_q({\mathfrak g})$ with $q$ a root of unity by the ideal of elements annihilated by representations with zero quantum dimension.  
 Since the ideal is not a coideal, the quotient algebra $A$ is not a bialgebra in the natural way. For ${\mathfrak sl}_2$, Mack and Schomerus defined a weak coproduct on $A$ and an associator derived from the fusion rules of the semisimple fusion category ${\mathcal C}({\mathfrak sl}_2, q)$,  a counit and a natural antipode  derived from $U_q({\mathfrak sl}_2)$, and have verified the axioms of a weak quasi-Hopf algebra on their example.   Their construction  is based on the special multiplicity free fusion rules of the category. In \cite{CGP}  we have extended Mack and Schomerus construction to all Lie types and to all roots of unity of the unitary type--i.e. the roots of unity appearing in the statement of our main Theorem \ref{new}--and we have derived    algebraic and analytic
  properties that were overlooked in the literature.

 \subsection{Weak Hopf algebras}    The main algebraic property that our quotient algebra $A_W({\mathfrak g}, q)$ of $U_q({\mathfrak g})$ satisfies is that of being the Hopf counterpart of a weak quasi-Hopf algebra, in a position completely analogous to Hopf algebras among quasi-Hopf algebras. For this reason, we shall refer to $A_W({\mathfrak g}, q)$ as a {\it weak Hopf algebra}.
 We have established the { weak Hopf algebra} structure of $A_W({\mathfrak g}, q)$ building upon   the quantum group fusion category ${\mathcal C}({\mathfrak g}, q)$ and Wenzl's linear fibre functor $W$ associated to it
  \cite{Wenzl}.

  Here we limit ourselves to recall the axioms of the definition of weak Hopf algebra which mostly differ from those of Hopf algebras.
 For the complete properties of antipode, counit and $R$-matrix of weak Hopf algebras we refer the reader to Sects. 9--14 of \cite{CGP}.

  \begin{defn}\label{defn_wh}
  Let $A$ be a semisimple unital complex algebra endowed with  a  coproduct $\Delta:A\to A\otimes A$ given by a not necessarily unital homomorphism and a counit satisfying the usual counit axioms for bialgebras.
The tensor product  $\rho\otimes\sigma$ of two f.d. representations     of $A$ is  the representation   defined by $\Delta$  restricted to the range
  of $\Delta(I)$ in   the tensor product space.
  
We call $A$ a {\it weak bialgebra} if the category ${\rm Rep}(A)$ of  representations of $A$   becomes a tensor category
  with associativity morphisms $$\alpha_{\rho, \sigma, \tau}: (\rho\otimes\sigma)\otimes\tau\to\rho\otimes(\sigma\otimes\tau),\quad \rho,\sigma,\tau\in{\rm Rep}(A)$$ defined by the action of $1\otimes \Delta(\Delta(I))\Delta\otimes 1(\Delta(I))$ between the representation spaces, with $I$ the unit of the algebra, and $1$ the identity map. In particular, if the coproduct is unital, $A$ is a coassociative bialgebra.
  
  A {\it weak Hopf algebra} is a weak bialgebra endowed with an antipode $S:A\to A$ satisfying the   axioms of a   Hopf algebra antipode.
  \end{defn}
  
  Among weak quasi-bialgebras, weak bialgebras admit a characterization in terms of weak tensoriality of their fibre functor, see Theorem 9.10 in \cite{CGP}.

\begin{prop} A semisimple    weak quasi-bialgebra $A$ is a weak bialgebra if and only if  the fibre functor ${\mathcal F}: {\rm Rep}(A)\to{\rm Vec}$ is
a {\it weak tensor functor}. By this we mean that    
the associated natural transformations  given by projection and inclusion, respectively:
$$F_{\rho, \sigma}: {\mathcal  F}(\rho)\otimes {\mathcal F}(\sigma)\to {\mathcal F}({\rho\otimes \sigma})$$ and $$G_{\rho, \sigma}: {\mathcal F}({\rho\otimes \sigma})\to {\mathcal F}(\rho)\otimes {\mathcal F}(\sigma)$$ 
  satisfy for all objects $\rho$, $\sigma$, $\tau$,

\begin{equation}\label{right_inverse}
F_{\rho, \sigma}\circ G_{\rho, \sigma}=1_{{\mathcal F}({\rho\otimes \sigma})}\end{equation}
\begin{equation}\label{naturality}
F_{\rho', \sigma'}\circ {\mathcal F}(S)\otimes{\mathcal F}(T)={\mathcal F}(S\otimes T)\circ F_{\rho,\sigma},\quad
{\mathcal F}(S)\otimes{\mathcal F}(T) \circ G_{\rho,\sigma}=G_{\rho',\sigma'}\circ {\mathcal F}(S\otimes T)
\end{equation}
for  objects $\rho$, $\sigma$, $\rho'$, $\sigma'\in{\mathcal C}$ and morphisms $S:\rho\to\rho'$, $T:\sigma\to\sigma'$.

\begin{equation}
{\mathcal F}(\alpha_{\rho,\sigma,\tau})=F_{\rho, \sigma\otimes\tau}
\circ 1_{{\mathcal F}(\rho)}\otimes F_{\sigma,\tau}\circ \alpha'_{{\mathcal F}(\rho), 
{\mathcal F}(\sigma), {\mathcal F}(\tau)}\circ G_{\rho, \sigma}\otimes 1_{{\mathcal F}(\tau)}\circ G_{\rho\otimes\sigma,\tau}\label{wt1}
\end{equation}
\begin{equation}{\mathcal F}(\alpha_{\rho,\sigma,\tau}^{-1})=  F_{\rho\otimes\sigma,\tau}\circ  F_{\rho,\sigma}\otimes 1_{{\mathcal F}(\tau)}\circ 
{\alpha'}^{-1}_{{\mathcal F}(\rho), 
{\mathcal F}(\sigma), {\mathcal F}(\tau)}\circ 1_{{\mathcal F}(\rho)}\otimes G_{\sigma, \tau} \circ G_{\rho, \sigma\otimes\tau},\label{wt2}\end{equation}
where $\alpha$ are the associativity morphisms of ${\rm Rep}(A)$ defined by the associator of $A$, and
$\alpha'$ are the (trivial) associativity morphisms of ${\rm Vec}$.

\end{prop}

   We refer to Sect. 6, 7, 10 in \cite{CGP} for the general definition of weak quasi-Hopf algebra, their $R$-matrix and how to obtain a rigid braided tensor category from it.
The structure of weak quasi-Hopf algebra with $R$-matrix is stable under a twist operation--first introduced by Drinfeld for quasi-Hopf algebras.
The twist of a weak quasi-Hopf algebra with $R$-matrix is automatically another weak quasi-Hopf algebra with $R$-matrix, and provides an equivalent rigid  braided tensor category.

\subsection{Unitary coboundary weak (quasi-)Hopf algebras}
 We next describe the properties concerning the unitary structure of weak quasi-Hopf algebras that are most relevant for this paper. The notion of a \textit{compatible coboundary unitary (or more generally Hermitian) weak quasi-Hopf algebra} is given in Def. 30.1 of \cite{CGP} as a strengthening of the more relaxed notion given in Def. 28.1 of the same paper. In the stronger case, a compatibility property between the opposed and adjoint structures of the weak quasi-bialgebra is required. In this paper, we focus on the stronger definition and tacitly assume the compatibility property. Regarding the weaker definition, its relevance relies on the categorical characterization extending the notion of a unitary symmetric functor for symmetric tensor $C^*$-categories to coboundary symmetries, by Theorem 29.1 in \cite{CGP}. On the other hand, the stronger definition allows for the abstract Drinfeld-Kohno Theorem 30.4 of \cite{CGP}

\begin{defn}\label{Hermitian_ribbon_wqh} A {\it (compatible) unitary coboundary weak (quasi-)Hopf   algebra} $A$ is defined by the following data:
\begin{itemize}
\item[{\rm        a)}]    A  unital weak (quasi-)Hopf  algebra $A$ endowed with a $C^*$-algebra
involution 
  with an antipode $S$  ($(S, \alpha, \beta)$ in the quasi-Hopf case)
\item[{\rm        b)}]   a
ribbon structure $(R, v)$ for $A$ associated to $(S, \alpha, \beta)$  
 such that the ribbon element $v\in A$
  is   unitary,  
\item[{\rm        c)}]     a unitary central square root $w\in A$
  of $v$ such that $\varepsilon(w)=1$, 
$S(w)=w$, 
\item[{\rm        d)}] $\tilde{A}=A^{{\rm        op}}$ as quasitriangular weak quasi-bialgebras,
\item[{\rm        e)}] $\overline{R}=Rw^{-1}\otimes w^{-1} \Delta(w)$ is positive.
 
  \end{itemize}
  
\end{defn}
  
 The element $R$ is the $R$-matrix; $v\in A$ is a central element defining the ribbon structure;
  the operator $\overline{R}$ is {\it Drinfeld coboundary matrix} and gives the $\Omega$-involution $\Omega=\overline{R}$ in the sense of Sect. 11, 12 in \cite{CGP}. The $\Omega$-involution induces a $^*$-involution on ${\rm Rep}(A)$ making it into a tensor $^*$-category. In the case where the $^*$-involution is induced by a coboundary matrix, the braided symmetry of ${\rm Rep}(A)$ induced by the $\overline{R}$-matrix is {\it unitary} by Theorem 28.9 in \cite{CGP}.
  
  The triple $(S, \alpha, \beta)$ describes the antipode of a weak quasi-Hopf algebra, and makes ${\rm Rep}(A)$ into a rigid category.
The presence of invertible elements $\alpha$, $\beta$ 
  is necessary in the quasi-cases, it arises from twist deformations of the usual antipode of a Hopf algebra.
  If $A$ is weak Hopf algebra, its antipode is required to satisfy the same axioms as a Hopf algebra antipode: is defined with $\alpha=\beta=I$; $\varepsilon$ denotes the counit and satisfies the same axioms as in the Hopf algebra case;   $A^{{\rm op}}$ and $\tilde{A}$,   denote the  opposite and adjoint   weak quasi-Hopf algebras. The opposite structure  $A^{{\rm        op}}$ is widely used in the literature, while
 the adjoint structure   $\tilde{A}$ is defined by coproduct   $\tilde{\Delta}(a)=\Delta(a^*)^*$, and similarly for the remaining structure, see Sect. 10 in \cite{CGP}.

Except for the $C^*$-structure of the algebra and the positivity requirement of the coboundary operator in e), the axioms
are satisfied by the $^*$-algebra of $U_q({\mathfrak g})$, which for this reason represents the motivating example of our definition.
These properties of $U_q({\mathfrak g})$ and positivity of a truncated coboundary operator operator of $U_q({\mathfrak g})$
  on certain special tensor products are  due to Wenzl \cite{Wenzl}, see also Sect. 31 in \cite{CGP}.
 
Compatible unitary weak quasi-Hopf  algebras in the sense of Definition \ref{Hermitian_ribbon_wqh} have been introduced and studied in Sect. 30 of \cite{CGP}.
They satisfy the stronger property d) than their weaker variant considered in Sect. 27 of the same paper, where d) is relaxed up to a so called {\it trivial twist}.  If the coproduct is unital, the two definitions coincide. We do not have examples in the weak case that distinguishes among the two definitions.

The weaker variant  is interesting for the fact that it admits a categorical characterization extending the notion of symmetric fibre functor
--a core notion in Doplicher-Roberts theory.--
to the coboundary symmetry induced by $\overline{R}$ in ${\rm Rep}(A)$. The stronger notion on the other hand allows  a formulation of an abstract Drinfeld-Kohno theorem, see Theorem 30.4 in \cite{CGP}, that we use in our main theorem.
In this paper we shall not need the weaker version, therefore we shall always   assume the compatibility condition given by property d) when not clearly stated.

In \cite{Wenzl} Wenzl has constructed a linear $^*$-functor 
$$W:{\mathcal C}({\mathfrak g}, q)\to {\rm Hilb}.$$
The following result has been shown in Theorem 2.5 in \cite{CGP}  (Theorem~9.2 of \cite{CP} in the type $A$ case).

\begin{thm}\label{wh}
The semisimple $^*$-algebra $A_W({\mathfrak g}, q):={\rm Nat}(W)$ of natural transformations of $W$ admits a natural structure
of  unitary compatible coboundary weak Hopf algebra derived from ${\mathcal C}({\mathfrak g}, q)$ and $W$.
\end{thm}

 Thanks to Theorem 30.4 in \cite{CGP}--our version of Drinfeld-Kohno theorem, the following result has been proved in Theorem 2.8 and Cor. 2.9 in \cite{CGP}

\begin{thm}\label{2.8} The Zhu algebra $A(V_{{\mathfrak g}_k})$ becomes a unitary coboundary weak quasi-Hopf algebras via a canonical Drinfeld twist isomorphism from $A_W({\mathfrak g}, q)$ for all Lie types. This structure induces a unitary modular fusion category structure on ${\rm Rep}(V_{{\mathfrak g}_k})$ equivalent to ${\mathcal C}({\mathfrak g}, q)$.
\end{thm}

\section{The main result} \label{5}\bigskip

Let ${\mathfrak g}$ be a complex simple Lie algebra and
let ${\rm Rep}(V_{{\mathfrak g}_k})$ be the module category of the corresponding affine vertex operator algebra $V_{{\mathfrak g}_k}$ at level $k\in{\mathbb N}$, endowed with the ribbon braided tensor structure introduced by Huang and Lepowsky.
Let $A(V_{{\mathfrak g}_k})$ be the Zhu algebra.

Let $V$ denote the fundamental representation of ${\mathfrak g}$ for each Lie type, defined in \cite{Wenzl}, (see also   Sect. \ref{12} 
for the classical types and $G_2$).
We also denote by $V$ the corresponding quantized representation of $U_q({\mathfrak g})$, $q=e^{i\pi/\ell d}$,
the associated object of the quantum group fusion category ${\mathcal C}({\mathfrak g}, q)$, 
  with $\ell=k+\check{h}$,   $\check{h}$ the dual Coxeter number of ${\mathfrak g}$.
  We use the same notation for 
    the corresponding
  representation of the vertex operator algebra $V_{{\mathfrak g}_k}$ and of the Zhu algebra  $A(V_{{\mathfrak g}_k})$.
  Similarly for the irreducible objects $V_\lambda$ of the involved   categories.

Let   $A_W({\mathfrak g}, q)$ denote the unitary compatible coboundary weak Hopf algebra associated to ${\mathcal C}({\mathfrak g}, q)$
via the construction of a weak tensor structure on Wenzl functor $W: {\mathcal C}({\mathfrak g}, q)\to{\rm Vec}$ 
 following  Theorem \ref{wh}.  
\bigskip

\begin{thm}\label{Finkelberg_HL}  
\begin{itemize}
We have that:
\item[(a)] 
The linear representation category ${\rm Rep}(A(V_{{\mathfrak g}_k}))$  
admits a natural structure of
ribbon braided tensor category ${\rm Rep}_{\rm HL}(A(V_{{\mathfrak g}_k}))$
obtained as a transport from  Huang-Lepowsky structure on ${\rm Rep}_{\rm HL}(V_{{\mathfrak g}_k})$ to ${\rm Rep}(A(V_{{\mathfrak g}_k}))$.
In this way,
Zhu's linear equivalence 
$$Z: {\rm Rep}_{\rm HL}(V_{{\mathfrak g}_k}) \to{\rm Rep}_{\rm HL}(A(V_{{\mathfrak g}_k}))$$ 
admits a natural structure of
ribbon braided tensor equivalence.

\item[(b)] 
The Zhu algebra $A(V_{{\mathfrak g}_k})$ becomes a twisted unitary compatible coboundary weak quasi-Hopf algebra with structure 
induced by  the   weak  Hopf algebra  $A_W({\mathfrak g}, q)$  
via Wenzl path and the construction of  the Drinfeld twist. Let
  ${\rm Rep}_{\rm QG}(A(V_{{\mathfrak g}_k}))$ be endowed with the corresponding 
  ribbon braided tensor structure. 
  
At the level of representation categories,
 ${\rm Rep}_{\rm QG}(A(V_{{\mathfrak g}_k}))$ becomes a unitary  coboundary modular tensor category in this way
   with the structure derived from ${\mathcal C}({\mathfrak g}, q)$ and the equivalence 
$${\rm Rep}_{\rm QG}(A(V_{{\mathfrak g}_k}))\xrightarrow{
}{\rm Rep}(A_W({\mathfrak g}, q)),$$
obtained in Theorem \ref{2.8}.
Moreover ${\rm Rep}_{\rm QG}(A(V_{{\mathfrak g}_k}))$ and ${\rm Rep}_{\rm HL}(A(V_{{\mathfrak g}_k}))$ have the same 
tensor product bifunctor.  
 
\item[(c)] 
The linear equivalence 
$$Z: {\rm Rep}_{\rm HL}(V_{{\mathfrak g}_k}) \to{\rm Rep}_{\rm QG}(A(V_{{\mathfrak g}_k}))$$ 
endowed with the same natural transformation as in (a) but w.r.t. the new ribbon braided tensor structure of  ${\rm Rep}_{\rm QG}(A(V_{{\mathfrak g}_k}))$ described in (b),
satisfies the ribbon braided tensor equivalence equations for   the ribbon  structure
for all objects,
the braid morphisms for pairs $(V_\lambda, V)$ and $(V, V_\lambda)$, and the associativity morphisms for triples
$$(V_\lambda, V, V), \quad (V, V_\lambda, V), \quad (V, V, V_\lambda)$$ with $V_\lambda$ an arbitrary irreducible object of the fusion category.

\item[(d)]  If  ${\mathfrak g}$
 is of one of the Lie types
   $A$, $B$, $C$, $D$, $G_2$     then 
     Zhu's equivalence in (c) is a ribbon braided tensor equivalence.   
     Specifically, the composition 
\begin{equation}\label{31.2}{\rm Rep}_{\rm HL}(V_{{\mathfrak g}_k}) \xrightarrow{Z}{\rm Rep}_{\rm QG}(A(V_{{\mathfrak g}_k}))\xrightarrow{}{\rm Rep}(A_W({\mathfrak g}, q))\xrightarrow{{\rm TK\  ribbon \ equiv}}{\mathcal C}({\mathfrak g}, q)\end{equation}

is a ribbon tensor equivalence by application of the indicated theorems.

\end{itemize}

\end{thm}

 We   notice that the equivalences described in (\ref{31.2}) are explicit.
 The   equivalence on the left is Zhu's linear equivalence endowed  with the natural transformations given by the restriction  to the zero modes  the braided tensor structure by Huang and Lepowsky.
 The fact that the braided tensor structure of the range of $Z$ induced by Huang-Lepowsky structure in (a) coincides with the structure  arising from quantum groups on the Zhu algebra,
 follows from  our work  in  (c), done in \cite{CGP}, and the first statement of part (d), established in this paper.
  The   equivalence on the right is the functor due to Wenzl enriched with our weak tensor structure by Theorem \ref{wh}.
    The middle equivalence   is an application of the explicit symmetry provided by our analogue of Drinfeld-Kohno theorem
30.4 in \cite{CGP}, and stated in parts (a) and (b) of Theorem \ref{2.8}. This symmetry is given by a Drinfeld twist explicitly defined by the action of the
  $R$-matrix, more precisely by the braided and ribbon structure in the two settings, that allow to see   in the vertex operator algebra setting
     same structure as that in the quantum group setting.
     
     \begin{rem}\label{Exceptional_types}
  The unitary modular tensor category structure on $\mathrm{Rep}_{\mathrm{QG}}(A(V_{\mathfrak{g}_k}))$ lifts to the entire category $\mathrm{Rep}(V_{\mathfrak{g}_k})$ via Zhu's linear equivalence $Z: \mathrm{Rep}(V_{\mathfrak{g}_k}) \to \mathrm{Rep}_{\mathrm{QG}}(A(V_{\mathfrak{g}_k}))$ and its quasi-inverse. This is established in the converse statement of Proposition \ref{transport} by virtue of Remark \ref{conversely}. Consequently, with respect to this transported structure, the Finkelberg--Kazhdan--Lusztig equivalence holds unconditionally across all Lie types.
      \end{rem}

\section{The strategy   of the  proof: unifying results from \cite{CGP} and the present paper}\label{6}
\bigskip

 In this section, we outline  the strategy of the proof of Theorem \ref{new}, and the more technical steps stated in Theorem \ref{Finkelberg_HL}, parts (a) through (d). The core of our strategy is to use the quantum gauge group $A_W(\mathfrak{g}, q)$ to transport the   categorical structures of quantum groups to the affine vertex operator algebra setting.

\subsection{Part (a): Transporting the Huang-Lepowsky ribbon braided tensor structure of ${\rm Rep}(V_{{\mathfrak g}_k})$ to the representation category  of the Zhu algebra $A(V_{{\mathfrak g}_k})$}

This first step  of Theorem 
\ref{Finkelberg_HL} follows from the application to Zhu's functor of    the  general construction. Let us consider  a 
(ribbon, braided) semisimple tensor category ${\mathcal C}$ together with a discrete complex algebra $A$ (a direct sum of full matrix algebras), a linear equivalence ${\mathcal C}\to{\rm Rep}(A)$ and an  inverse linear equivalence
${\rm Rep}(A)\to{\mathcal C}$. We  use these equivalences to transfer all the structure of  ${\mathcal C}$ to ${\rm Rep}(A)$, and make it into a (ribbon,  braided) tensor category in such a way that the given equivalences become (ribbon, braided) tensor equivalences. 

The definition of the ribbon braided tensor structure on ${\rm Rep}(A)$ is   given in Proposition \ref{transport} and verification of the axioms of (ribbon, braided) tensor category and the natural transformations on the equivalences and the needed properties   is straightforward. 
We have used a similar construction, in the reverse direction, and for the purpose of construction of unitary structures on tensor categories from unitary weak quasi-Hopf algebras,
  in Theorems 15.6 and 23.10 of \cite{CGP}. Proposition~\ref{transport} itself, and the remarks following it,  are stated  here in the form the present argument requires.

The  application to   ${\rm Rep}(A(V_{{\mathfrak g}_k}))$, made into  a ribbon braided tensor category from Huang-Lepowsky structure
on ${\rm Rep}(V_{{\mathfrak g}_k})$,   corresponds to passing from the structure on a category of infinite dimensional modules to the category of their zero modes. We have explicitly described this construction  in Sect. 37 in \cite{CGP}. 
This construction  
 was    also considered   in \cite{McRae}.
 
\textbf{Main result of part (a):} Through this identification, we establish that the representation category of the Zhu algebra $A(V_{\mathfrak{g}_k})$--which determines the initial terms of the primary fields--inherits the structure of a ribbon braided tensor category.

\subsection{Part (b): Making ${\rm Rep}(A(V_{{\mathfrak g}_k})$ into a unitary ribbon braided tensor category via a Drinfeld twist method from
$A({\mathfrak g}, q)$}

This step is independent of part (a) and makes ${\rm Rep}(A(V_{{\mathfrak g}_k}))$ into a unitary ribbon braided tensor category equivalent to the quantum group fusion category ${\mathcal C}({\mathfrak g}, \ell)$ for all Lie types.
The core strategy   consists in   developing  an abstract analogue of the Drinfeld-Kohno theorem (see Theorem 30.4 in \cite{CGP}). We apply our Drinfeld twist, composed with a quantization isomorphism derived from Wenzl's work \cite{Wenzl} to establish a bridge between the unitary coboundary weak Hopf algebra $A_W(\mathfrak{g}, q)$ and the Zhu algebra $A(V_{{\mathfrak g}_k})$. 
The method was originally used by Drinfeld for Drinfeld-Kohno theorem
\cite{Drinfeld_cocommutative}, \cite{Drinfeld_quasi_hopf}, \cite{Drinfeld_galois}.  

 Specifically, this step splits in various parts. The first main part is the construction of the class of finite dimensional {\it unitary coboundary weak Hopf $C^*$-algebras}  $A_W({\mathfrak g}, q)$ recalled in Theorem \ref{wh} for all the
 Lie types and positive integer levels. This construction extends with different methods the previous examples of \cite{CP} in the type $A$ case.
 The unitary coboundary structure of  $A_W({\mathfrak g}, q)$ naturally arises from the unitary structure of the fusion category of quantum groups at roots of unity, associated to the specific linear fibre
 functor $W: {\mathcal C}({\mathfrak g}, q)\to{\rm Vec}$,  from the work of Wenzl \cite{Wenzl}.
 
 The second part of the proof for (b) is the construction of the Drinfeld twist $F$ defined as a suitable square root
 of $\overline{R}$.   Therefore in our case, the unitary structure of the fusion category of quantum groups at roots of unity in \cite{Wenzl} is used
 in the construction of the Drinfeld twist of   $A_W({\mathfrak g}, q)$. The third part is  the development of
the mantioned abstract analogue of Drinfeld-Kohno theorem 30.4 in \cite{CGP}, which we apply to our twist composed with a quantization isomorphism 
described in \cite{Wenzl} for the unitary structure of irreducible representations of  ${\mathcal C}({\mathfrak g}, q)$.  

 \textbf{Main result of part (b):} 
 This identifies the twist of
  an isomorphic image of the $C^*$-algebra $A_W({\mathfrak g}, q)$ with a corresponding quotient of $U({\mathfrak g})$ that in turn identifies with the Zhu algebra
 $A(V_{{\mathfrak g}_k})$, and in this way becomes a {\it unitary coboundary weak quasi-Hopf algebra} with associator given by a $3$-coboundary. 

  Moreover, our Drinfeld-Kohno Theorem, similarly to the original Drinfeld theorem,
 explicitly turns the braided symmetry of $A_W({\mathfrak g}, q)$ into the braided symmetry
 known in the case of affine Lie algebras and affine vertex operator algebras, at least on all those fusion tensor products of representations
 of the form $V\otimes V_\lambda$, with $V$ the fundamental and $V_\lambda$ any irreducible. This represents the first indication that our twisted ribbon braided structure corresponds to Huang-Lepowsky ribbon braided tensor structure, and we shall build on it for the proof of  the 
ribbon braided tensor equivalence.

 This induces the structure of a unitary
 ribbon braided tensor category
 on the module category of   $V_{{\mathfrak g}_k}$, via transport from   the Zhu algebra $A(V_{{\mathfrak g}_k})$ following part  (a) in the reverse direction.
 In particular,  we get an independent construction of the associativity morphisms and proof of rigidity of ${\rm Rep}(V_{{\mathfrak g}_k})$, directly arising from the quantum group fusion category via the equivalence induced by the twist.

 By construction, this structure is ribbon braided tensor equivalent to
 the structure arising from the quantum group fusion category.

 \subsection{Part (c): Constructing the ribbon braided tensor equivalence between ${\mathcal C}({\mathfrak g}, q)$ and ${\rm Rep}(V_{{\mathfrak g}_k})$}
 
This third step initiates the comparison between Huang-Lepowsky braided tensor structure on ${\rm Rep}_{\rm HL}(V_{{\mathfrak g}_k})$ and our structure on ${\rm Rep}_{\rm QG} (A(V_{{\mathfrak g}_k}))$. 
The axioms of braided tensor equivalence are verified in \cite{CGP} by a direct comparison  for the special braiding morphisms and   the special associativity morphisms listed in (c) on the Zhu algebra, passing through (a). Up to this point we do not need any knowledge on
  the space of intertwiners in the two categories, the linear equivalence is an isomorphism between any two corresponding morphism spaces by the Drinfeld twist construction.
  
  The necessary general theory
of semisimple weak quasi-Hopf algebras, weak  Hopf algebras, and the theory of Drinfeld twist treated in
 Sections  4--7, and  9--17 in \cite{CGP}. 
  This  abstract theory splits in two parts. Theory on which   part (b), (c) of Theorem \ref{Finkelberg_HL} is based
 is developed in Sections 4--7, and 9--16. 
 The applicative part of the proof of part (b), (c) is
   developed in  Sections 19, 20, 22--38 of that paper.

   \subsection{Part (d): Completing the proof of the braided tensor equivalence for the classical Lie types and $G_2$}
   
 The abstract aspects that lead to the proof of part (d) of Theorem \ref{Finkelberg_HL}  are the main subject of  Section 8 of \cite{CGP}.
 In Sect. 8 of \cite{CGP} we state two main abstract uniqueness results for   associativity morphisms and braided symmetries in a semisimple pre-tensor category with a generating object. These results, along with their complete proofs, are presented below as Theorem \ref{claim0} and its extension Theorem  \ref{claim1}.   In the subsequent sections we explain how   Theorem   \ref{claim1} applies to establish a first step to part (d) of Theorem \ref{Finkelberg_HL}. 
These are the main   abstract results to obtain the complete braided tensor equivalence, and rely on  key aspects of the discovery of quantum groups. 
In Sect. \ref{13} of this paper we apply these theorems to complete the proof of part (d) of Theorem \ref{Finkelberg_HL}. In this case we work on the comparison between the two structures transported to the side of the quantum groups, $A_W({\mathfrak g}, q)$, and we need knowledge on the generating property by the braid group between tensor powers of the fundamental representation, to apply our uniqueness theorem.
 \bigskip

\section{Outline of the  Tannakian  construction of $A_W({\mathfrak g}, q)$
 }\label{7}\bigskip

\noindent Let  $({\mathcal C}, \otimes, \iota, \alpha)$ be a semisimple  complex   tensor category   
  (possibly   
braided,  rigid, ribbon)  and let
$${\mathcal F}: {\mathcal C}\to{\rm Vec}$$ be a faithful
linear functor. 
There are  two natural ways to make the representation category ${\rm Rep}(A)$ of $A={\rm Nat}_0({\mathcal F})$,   the algebra of natural transformations of ${\mathcal F}$ with finite support,
into a tensor category   with the same type of structure as ${\mathcal C}$,  equivalent to ${\mathcal C}$ with all the   structure.

The first way is traditionally referred to as 
Tannaka-Krein duality. We refer to the second way   as a {\it transport} of the structure of ${\mathcal C}$ and will be outlined in the next section.

To implement Tannaka-Krein duality, one needs to construct
  a   weak, quasi-tensor structure $(F, G)$ on ${\mathcal F}$ (i.e. $F$ and $G$ satisfy (\ref{right_inverse})  and then apply
   Theorem 7.6 in \cite{CGP} to obtain    a
weak, quasi-bialgebra structure on $A={\rm Nat}_0({\mathcal F})$. Coproduct $\Delta$,  counit $\varepsilon$,  associator  $\Phi$ are respectively defined
for any $a \in A$ by:
$$ \Delta(a)_{\rho,\sigma} = G_{\rho,\sigma} a_{\rho \otimes \sigma} F_{\rho,\sigma}, \quad \varepsilon(a) = f^{-1} a_\iota f,$$
where $f: \mathbb{C} \to {\mathcal F}(\iota)$,
\begin{equation}\label{associativity}\Phi_{\rho,\sigma, \tau}=1_{{\mathcal F}(\rho)}\otimes G_{\sigma, \tau}\circ G_{\rho, \sigma\otimes\tau}\circ{\mathcal F}(\alpha_{\rho, \sigma,\tau})\circ F_{\rho\otimes\sigma,\tau}\circ F_{\rho, \sigma}\otimes 1_{{\mathcal F}(\tau)}.
\end{equation}
Antipode, $R$-matrix, ribbon element of $A$ follow respectively from rigidity, braided, ribbon structure of ${\mathcal C}$.
The coproduct, associator, $R$-matrix, and ribbon element lie in corresponding multiplier algebras of $A$, or its tensor powers.

Regarding the associator, Tannaka-Krein duality is particularly simple in cases where ${\mathcal C}$ is  {\it strict}: $\alpha=1$.
In this case the associator $\Phi$   depends only on $F$ and $G$.  
      
Fusion categories that arise from quantum groups at roots of unity ${\mathcal C}({\mathfrak g}, q)$ can be realized
   as   strict tensor categories. To the author's knowledge, the argument of this result is due to the combined methods of 
   \cite{Gelfand_Kazhdan} and \cite{Wenzl}. For an explicit proof see Theorem 5.4 in \cite{CP}.  
   A natural  linear functor
   $$W: {\mathcal C}({\mathfrak g}, q)\to{\rm Hilb}$$
     has been naturally constructed in \cite{Wenzl} for the roots of unity of interest for the aims of \cite{CGP}, i.e. 
   those that are related to the affine Lie algebras at positive integer levels.
    In this example, a main point of our work is the construction of $F$ and $G$   as fusion projections and inclusions on and from subspaces.
    
    It follows that the associator of $A_W({\mathfrak g}, q)$ by Tannakian duality becomes that of a weak Hopf algebra
    in the sense of definition \ref{defn_wh}.
By Theorem \ref{wh}, $A_W({\mathfrak g}, q)$ can be constructed with the structure of a compatible unitary coboundary weak Hopf $C^*$-algebra.
    \bigskip

   \section{The induced Huang-Lepowsky structure on 
   ${\rm Rep}(A(V_{{\mathfrak g}_k}))$ (Proof of Th. \ref{Finkelberg_HL} (a))}\label{8}

  A second constructive method to endow the representation category 
${\rm Rep}(A)$ with a ribbon braided tensor category structure from a given semisimple ribbon braided tensor category 
$({\mathcal C}, \otimes, \iota, \alpha, c, v)$ is applicable when ${\mathcal C}$ is equipped with a faithful linear functor
$${\mathcal F}: {\mathcal C}\to{\rm Vec}.$$
In this case, one first considers the discrete algebra $A={\rm Nat}_0({\mathcal F})$ and then transfers
 the structure of ${\mathcal C}$ directly to ${\rm Rep}(A)$.

This transfer does not require a weak quasi-tensor structure on ${\mathcal F}$. Instead, it is constructed using the linear equivalence ${\mathcal E}: {\mathcal C}\to {\rm Rep}(A)$ associated to ${\mathcal F}$ and an inverse linear equivalence ${\mathcal S}: {\rm Rep}(A)\to {\mathcal C}$. 

\begin{prop}\label{transport} Let ${\mathcal E}: {\mathcal C}\to {\rm Rep}(A)$ be a  linear equivalence with the representation category of a discrete algebra $A$, and let 
 ${\mathcal S}: {\rm Rep}(A)\to {\mathcal C}$ be an inverse linear equivalence.  Then
the braided tensor category structure on ${\rm Rep}(A)$ is explicitly defined by transporting the structure from ${\mathcal C}$: for any objects $U, V, W$ in ${\rm Rep}(A)$, a tensor product object is given by
\begin{equation}
U \otimes_A V := {\mathcal E}({\mathcal S}(U) \otimes_{\mathcal C} {\mathcal S}(V)),
\end{equation}
a tensor product morphism by a similar formula, and the associativity, braiding, ribbon morphisms are respectively defined by
\begin{equation}
\alpha^A_{U,V,W} := {\mathcal E}(1_{\mathcal{S}(U)}\otimes\eta^{-1}_{\mathcal{S}(V)\otimes \mathcal{S} (W)}  \circ \alpha^{\mathcal C}_{{\mathcal S}(U), {\mathcal S}(V), {\mathcal S}(W)} \circ 
\eta_{\mathcal{S}(U)\otimes \mathcal{S}(V)}\otimes 1_{\mathcal{S}(W)}), 
\end{equation}
where $\eta:{\mathcal S}{\mathcal E}\to 1$ is an invertible natural transformation,
 \begin{equation} c^A_{U,V} := {\mathcal E}(c^{\mathcal C}_{{\mathcal S}(U), {\mathcal S}(V)}), \quad v^A_U:={\mathcal E}(v^{\mathcal C}_{{\mathcal S}(U)}).
\end{equation}

With this structure, the pair $({\mathcal S}, S)$, where $S_{U, V}=\eta^{-1}_{{\mathcal S}(U)\otimes {\mathcal S}(V)}$, constitutes an equivalence of ribbon braided tensor categories.

\end{prop}

\begin{proof}
The proof follows from straightforward verification.
\end{proof}

  Fundamental examples of such a functor ${\mathcal F}$ occur in conformal field theories (e.g., in the setting of models of conformal nets and vertex operator algebras) as the {\it minimum energy functor}. 
   In the case of vertex operator algebras, we refer to the braided tensor category structure of modules of vertex operator algebras by Huang and Lepowsky. For conformal nets we refer  to the unitary braided tensor category structure of localized endomorphisms of Doplicher-Haag-Roberts.

  \subsection{Zhu's linear equivalence ${\mathcal E}: {\rm Rep}(V)\to {\rm Rep}(A(V))$ and its   inverse ${\mathcal S}$.}
In the setting of vertex operator algebras, the minimum energy functor is also known as {\it Zhu's functor} \cite{Zhu}. This
functor which assigns to a $V$-module $W$ its lowest energy subspace $W(0)$. Zhu explicitly constructed an associative   algebra $A(V)$,  directly from the modules of the vertex operator algebra $V$. This is the well-known Zhu algebra. 

Furthermore, under suitable semisimplicity conditions, Zhu defined a linear equivalence ${\mathcal E}: {\rm Rep}(V)\to {\rm Rep}(A(V))$ between module categories of $V$ and of $A(V)$, acting as $W \mapsto W(0)$, alongside an explicit natural inverse equivalence ${\mathcal S}: {\rm Rep}(A(V))\to {\rm Rep}(V)$. 
Thus $A(V)$ is   isomorphic to ${\rm Nat}_0({\mathcal F})$.
The infinite dimensional $V$-module ${\mathcal S}(U)$ associated to a module of the Zhu algebra is called {\it Zhu's induced module.}

  Thus, in the case of vertex operator algebras satisfying semisimplicity conditions, the analysis of the category of infinite-dimensional representations of the vertex operator algebra $V$ is constructively reduced to the category of finite-dimensional representations of the Zhu algebra $A(V)$. For affine vertex operator algebras at a positive integer level $k$, this explicit transfer constitutes the passage described in part (a) of the statement of Theorem \ref{Finkelberg_HL}.

 \begin{rem} It is important to recall that  Zhu's inverse linear equivalence 
  ${\mathcal S}: {\rm Rep}(A(V))\to
  {\rm Rep}(V)$
   is a {\it right inverse} of ${\mathcal E}$: ${\mathcal E}\circ {\mathcal S}=1$, that is the minimum energy subspace of ${\mathcal S}(U)$ is precisely $U$. This fact implies that for an affine vertex operator algebra $V_{{\mathfrak g}_k}$ the tensor product of   $A(V_{{\mathfrak g}_k})$-modules derived in Prop. \ref{transport} from Huang-Lepowsky tensor structure 
   ${\rm Rep}_{\rm HL}(V_{{\mathfrak g}_k})$ can be identified as a
  {\it canonical addendum submodule} of the simple Lie algebra ${\mathfrak g}$  for all Lie types. The   realization of the tensor structure of the Zhu algebra modules derived from Huang-Lepowsky tensor structure described in the following remark, is at the core of our proof of rigidity
  of ${\rm Rep}_{\rm HL}(V_{{\mathfrak g}_k})$. Our proof of rigidity of    ${\rm Rep}_{\rm HL}(V_{{\mathfrak g}_k})$
  arises from   the antipode of the Hopf algebra $U_q({\mathfrak g})$ at roots of unity and the tensor ideal nature of the subcategory of negligible modules. Indeed, 
  building on Wenzl work in the framework of quantum group fusion categories   \cite{Wenzl},  we   eventually show that the same transferred structure is obtained from a Drinfeld twist method from the unitary   coboundary weak Hopf algebra $A_W({\mathfrak g}, q)$ associated to the quantum group fusion category.
  \end{rem}

  \begin{ex}\label{special} {\it The   tensor product structure
  of ${\rm Rep}(A(V_{{\mathfrak g}_k}))$ induced by Huang-Lepowsky tensor structure and Zhu's linear equivalence.} Let us consider ${\mathcal C}={\rm Rep}(V_{{\mathfrak g}_k})$, the module category of an affine vertex operator algebra at positive integer level $k$.
   Since
   Zhu's inverse linear equivalence  ${\mathcal S}$  is a right inverse of ${\mathcal E}$, the minimum energy subspace of Zhu's induced module ${\mathcal S}(U)$ is precisely $U$. It follows that if $V_{\rm fund}$ is the   representation of the Zhu algebra $A(V_{{\mathfrak g}_k})$--identified with a quotient of the universal enveloping algebra $U({\mathfrak g})$--corresponding to the fundamental representation of   ${\mathfrak g}$   and $U$ 
  is any other irreducible representation of $A(V_{{\mathfrak g}_k})$,  then the minimum energy subspace of the Huang-Lepowsky tensor product module
  ${\mathcal S}(U)\otimes_{\rm HL}{\mathcal S}(V_{\rm fund})$ identifies with  a canonical addendum submodule of the usual complex tensor product ${\mathfrak g}$-module $U\otimes_{\mathbb C} V_{\rm fund}$   identified by the range   of a ${\mathfrak g}$-intertwining idempotent $F_0$ with dominant weights in the Weyl alcove of level $k$. This can be done by paralleling to the classical case the analogous  arguments of \cite{Wenzl} on the analysis of fusion of the fundamental representation in the quantum case     for 
  ${\mathfrak g}\neq E_8$.
   This submodule is moreover a direct sum of isotypic   ${\mathfrak g}$-representations. Indeed, the minimum energy 
  subspace of  ${\mathcal S}(U)\otimes_{\rm HL}{\mathcal S}(V_{\rm fund})$
   identifies with the $A(V_{{\mathfrak g}_k})$-module  defined by the   structure map $F_0$  obtained directly from the Huang-Lepowsky tensor product structure map $F$ acting between   Zhu's induced
    infinite-dimensional $V_{{\mathfrak g}_k}$-modules. Consequently,
by definition, the
   tensor product  $U\otimes_A V_{\rm fund}$ is  precisely this ${\mathfrak g}$-submodule of $U\otimes_{\mathbb C} V_{\rm fund}$ defined by $F_0$. These special tensor products play the role of {\it building blocks} of the braided tensor structure.
   Fundamental representations for all Lie types and their properties are discussed in Sections \ref{9}, \ref{12}.
   
  \end{ex}

\begin{rem}{\it The $E_8$-case.} The construction of  $F_0$ as an idempotent on the building block tensor products of ${\mathfrak g}$-modules of Example \ref{special} is linked to an analogous construction in \cite{Wenzl} of an idempotent
onto a submodule of a corresponding  tensor product module of   $U_q({\mathfrak g})$ at the   root of unity corresponding to the level $k$,
via his de-quantization continuous curve. This is necessary especially in the
$E_8$-case because the algebraic isotypic decomposition fails.
\end{rem}

  To prove Theorem \ref{Finkelberg_HL}, we first 
construct an identification between a quantum group fusion category 
at a fixed root of unity and the corresponding affine VOA fusion category 
as pre-tensor categories in the sense introduced in \cite{CGP} (i.e. as categories endowed only with a tensor product bi-functor and a tensor unit). We then identify their full structures on the building block tensor products, and eventually to the whole category.

As shown in   \cite{CGP} starting from Section 34, 
the Huang and Lepowsky theory, combined with previous work by Frenkel and Zhu, Wenzl's work, and our Drinfeld twist method,
naturally provides ${\mathcal F}$ with a specific weak quasi-tensor structure
 $(F_0, G_0)$, where $F_0$ is again explicitly constructed as the $0$-grade part of the 
 Huang-Lepowsky tensor structure map $F$, and extends the map $F_0$ described in Example
 \ref{special} from the special pairs  to all pairs of modules.

\begin{rem}  Specifically, it easily follows from Theorem 34.7 of \cite{CGP} that the above ${\mathfrak g}$-submodule property extends from building block tensor products to arbitrary tensor products $U\otimes_A W$ 
  of $A(V_{{\mathfrak g}_k})$-modules: the tensor product module
    $U\otimes_{A(V_{{\mathfrak g}_k})} W$  associated to Huang-Lepowsky tensor structure of ${\rm Rep}(V_{{\mathfrak g}_k})$
    following Proposition \ref{transport} can be identified with the range of an idempotent
   $(F_0)_{U, W}$ on the ${\mathfrak g}$-module  $U\otimes_{\mathbb C} W$ onto an invariant submodule.
In this way we construct natural weak quasi-tensor structures $({\mathcal F}, F_0, G_0)$, with $G_0$ the inclusion map, 
associated with an explicit isomorphic twist of the weak Hopf algebra $A_W({\mathfrak g}, q)$. Conversely, starting from the Zhu algebra $A(V_{{\mathfrak g}_k})$ as endowed with the   structure of a unitary coboundary weak quasi-Hopf algebra twisted isomorphic to the unitary compatible weak Hopf algebra $A_W({\mathfrak g}, q)$, these addendum ${\mathfrak g}$-submodules $U\otimes_{A(V_{{\mathfrak g}_k})} W$ define a self-contained   structure of a unitary ribbon rigid tensor category on ${\rm Rep}(A(V_{{\mathfrak g}_k}))$--equivalent to ${\rm Rep}(A_W({\mathfrak g}, q))$--that coincides with the structure induced by Huang and Lepowsky on tensor products. We   eventually show in Section \ref{12} that this identification extends to the whole category and the whole structure as made precise in Theorem \ref{Finkelberg_HL}.
\end{rem}

By applying a cohomological uniqueness result for braided tensor structures on any category 
with a generating representation (Theorem \ref{claim1}), combined 
with braid group duality properties established for quantum groups, we uniquely extend
 this identification to all variables. A more detailed overview of
  this constructive proof is provided in the following sections.

  Because  ${\mathcal C}({\mathfrak g}, q)$ is a unitary, rigid, modular tensor 
category, we   obtain by transport  unitarity, rigidity, and modularity of ${\rm Rep}(V_{{\mathfrak g}_k})$
 directly from the quantum group fusion category for all Lie types. 
 In particular we obtain an independent construction of rigid modular tensor structure on ${\rm Rep}(V_{{\mathfrak g}_k})$. By part (d) of Theorem \ref{Finkelberg_HL}, this in particular gives a new proof of
 Huang's rigidity and modularity   \cite{Huang2}
 for affine vertex operator algebras at a positive integer level for the stated Lie types and Huang-Lepowsky tensor structure derived from the quantum group fusion category. This is stated more precisely in the following remark.
 
        \begin{rem}\label{conversely}
 Conversely, assume that $\mathrm{Rep}(A)$ is equipped with the structure of a ribbon braided tensor category, and let $\mathcal{E}: \mathcal{C} \to \mathrm{Rep}(A)$ be a linear equivalence of semisimple linear categories with an inverse equivalence $\mathcal{S}$. By interchanging the roles of $\mathcal{C}$ and $\mathrm{Rep}(A)$ in Proposition \ref{transport}, $\mathcal{C}$   inherits a ribbon braided tensor category structure. 
 
 Assume moreover that $A$ is a discrete unitary weak quasi-Hopf $C^*$-algebra in the sense of Def. 11.5 in \cite{CGP} with $R$-matrix and unitary ribbon element.
 Then   $\mathrm{Rep}(A)$ becomes a  unitary ribbon braided tensor category  under a natural compatibility condition
 between the $\Omega$-involution of $A$ and the $R$-matrix, by Prop. 13.2 in \cite{CGP}.
 If we know that  ${\mathcal C}$ has the structure of a linear $C^*$-category then ${\mathcal C}$ upgrades to a unitary ribbon braided tensor category, by Theorem 15.7 of \cite{CGP}.
 This reverse transport of structure is employed  to construct the unitary ribbon braided tensor category structure on $\mathrm{Rep}(V_{\mathfrak{g}_k})$ from a given linear $C^*$-structure on this category and  the corresponding structure on the module category of the Zhu algebra $A(V_{\mathfrak{g}_k})$, which is in turn derived via the Drinfeld twist method in Sections \ref{9} and \ref{10}. This corresponds precisely to the structure mentioned in the first part of the introduction, originating from the quantum group fusion category $\mathcal{C}(\mathfrak{g}, q)$.   
 \end{rem}

\section{Unitarity transport and   the weak quasi-Hopf structure of the Zhu algebra}\label{9}

  A central method for establishing parts (b), (c), and (d) of Theorem \ref{Finkelberg_HL} relies on the structural properties of the fundamental generating representation $V_{\rm fund}$ across all Lie types, precisely the building block analyzed in Example \ref{special}. For (b), (c), we manage the quantum group side of the equivalence by utilizing the representation theory of $U_q({\mathfrak g})$ at roots of unity developed by Wenzl \cite{Wenzl}. 
  
  To bridge the non-semisimple setting of $U_q({\mathfrak g})$ with the semisimple setting of affine vertex operator algebras, our strategy is to verify the exact coincidence of the braided tensor structures on specific generating variables, and then algebraically extend this coincidence to the entire category.

Specifically, parts (b), (c) of Theorem \ref{Finkelberg_HL}  give the compatibility (we verify exact coincidence using the equivalence) of the two braided symmetries and associativity morphisms    all Lie types on the stated specific collection  of their variables,   in the two settings of    quantum groups and    vertex operator algebras. These collections
do not generate additively all the variables of the braiding morphisms or associativity morphisms, because they do not contain all the irreducible objects. 

The virtue of these special objects is that compatibility of the associativity morphisms and braiding morphisms on them is  verified in a direct way, 
based on  the same   structural properties of fusion of two
    corresponding generating representations in the two cases. In the quantum group case, fusion arises in the      setting of 
    non-semisimple representations of $U_q({\mathfrak g})$ considered in \cite{Wenzl}, while in the vertex operator algebra
    case fusion arises in a semisimple setting.

     To compare the braiding and associativity morphisms on the special variables, we use a Drinfeld twist that we construct
    trivialising, in a sense, the non-trivial unitary structure in the setting of representations of quantum groups at roots of unity to that of the affine vertex operator algebra. This idea emerges in the construction and analysis of the unitary structure of the fusion category of quantum groups at roots of unity of \cite{Wenzl}.
    
    The twist takes one unitary structure of a special tensor product arising in the quantum group setting from
     \cite{Wenzl}
    to the classical unitary structure corresponding to the compact real form of ${\mathfrak g}$
    (and this in turn remarkably leads to the construction of a unitary tensor  structure on the module category of the vertex operator algebra $V_{{\mathfrak g}_k}$
    in a uniform way).
This implies that  the twist leaves 
    the corresponding structural properties invariant on these special variables.   
 The following is a more detailed  description.

\subsection{The Drinfeld twist and the unitary coboundary weak quasi-Hopf structure on Zhu algebra}
Parts (b) and (c) of Theorem \ref{Finkelberg_HL} establish the exact coincidence of the two braided tensor structures on a distinguished collection of generating objects involving $V_{\rm fund}$. 

To compare the associativity and braiding morphisms rigorously, we must map the non-trivial unitary structure of the quantum group tensor product to the classical unitary structure corresponding to the compact real form of ${\mathfrak g}$, making in this way
${\rm Rep}(V_{{\mathfrak g}_k}$ into a {\it unitary} modular tensor category.

As detailed in Sections 31 and 32 of \cite{CGP}, the weak tensor structure of the Wenzl functor enables the construction of the unitary coboundary weak Hopf algebra $A_W({\mathfrak g}, q, \ell)$ directly from $U_q({\mathfrak g})$. To bridge this to the Zhu algebra representation category, we utilize Wenzl's continuous curve of $C^*$-tensor categories \cite{Wenzl} (see also Sect. 31 in \cite{CGP}). Parallel transport along this curve induces a natural analytic isomorphism of vector spaces:
\begin{equation}
\psi: V_\lambda(q) \longrightarrow V_\lambda(1).
\end{equation}
for every irreducible object $V_\lambda(q)$ of the quantum group fusion category, assigning an isomorphism with an irreducible object of the Lie algebra using Lusztig-Kashiwara basis \cite{Wenzl}. This isomorphism specifically identifies the representation modules $V(q)$ of the quantum gauge group $A_W({\mathfrak g}, q, \ell)$ with the minimum energy spaces $V(1)$ of the VOA (the modules of the Zhu algebra), providing a natural transformation between their respective fiber functors.

To transport the braided symmetry, we consider the Drinfeld coboundary matrix $\bar{R}$ derived from the quantum $R$-matrix of $A_W({\mathfrak g}, q, \ell)$. Following \cite[Eq. 29.17]{CGP}, $\bar{R}$ is defined by modifying the braiding with a square root $w$ of the ribbon element $v$:
\begin{equation}
\bar{R} = R\Delta(w)(w \otimes w)^{-1}.
\end{equation}
A crucial property established for the compatible unitary coboundary weak quasi-Hopf algebra $A_W({\mathfrak g}, q, \ell)$ is that this coboundary matrix $\bar{R}$ is self-adjoint by Theorem 28.9 in \cite{CGP}. Applying Theorem 30.4 in \cite{CGP}, an analogue of the original Drinfeld-Kohno theorem,
we construct a specific unitary Drinfeld twist $T$ as a square root of $\bar{R}$. 

Combining the analytic isomorphism $\psi$ with the twist $T$, we explicitly define the induced weak quasi-tensor structure on the Zhu algebra. A classical coproduct $\Delta_{0}$ of the Zhu algebra---which strictly agrees with the Huang-Lepowsky tensor product via the map $F_0$ of Example \ref{special}---is related to the quantum coproduct $\Delta$ via:
\begin{equation}
\Delta_0(a) = (\psi \otimes \psi) T\Delta(\psi^{-1}(a))T^{-1} (\psi \otimes \psi)^{-1} , \quad \forall a \in A(V_{{\mathfrak g}_k}).
\end{equation}

A similar twisted formula holds between the two $R$-matrices. The correct transformation of the squared $R$ matrix $R_{21}R$  first motivated the Drinfeld-Kohno theorem 30.4 of \cite{CGP}, similarly to the original Drinfeld-Kohno theorem. In this way, the $2$-coboundary defined by the ribbon elements
become compatible via the twist.
The framework of weak quasi-Hopf algebras thus allows to transport   the trivial weak Hopf algebra associator 
$$ \Phi_{\rm QG}=1\otimes\Delta(\Delta(I))\Delta\otimes 1(\Delta(I))$$ of $A_W({\mathfrak g}, q)$
to a $3$-coboundary associator of the Zhu algebra $A(V_{{\mathfrak g}_k})$. The twist $T$ acts as a 2-cocycle on the tensor product of the generating representations. The twisted associator $ \Phi_{\rm Z}$ (representing the Zhu algebra associator up to the isomorphism 
$\psi$) is obtained from the associator $\alpha$ of $A_W({\mathfrak g}, q)$ via the standard cohomological formula, explicitly revealing its 3-coboundary nature:
\begin{equation}\label{eq:twisted_associator}
\Phi_{\rm Z}= {1} \otimes T({1} \otimes \Delta)(T) \circ \Phi_{\rm QG} \circ (\Delta \otimes {1})(T^{-1})T^{-1} \otimes \mathbf{1}.
\end{equation}

\subsection{The equivalence of braiding and associator on the special variables} 
Specifically, the   coincidence of the associativity morphisms  and braided symmetry on special objects
 are verified   on the side
 of vertex operator algebras. 
  To see this,     we compare with the work by Frenkel and Zhu \cite{Frenkel_Zhu}. 
 We apply   Wenzl isomorphism and    Drinfeld twist construction and then we use
  the structural properties of the Huang-Leposwky associativity morphisms in Theorem 33.7 in \cite{CGP}. For the special   braiding
morphisms we take into account the formulas given by Toledano-Laredo in Theorem 33.11 in \cite{CGP}.  For the case of affine vertex operator algebras, see Remark    38.16 in the same paper.

  \subsection{Cohomological reduction via the pentagon and hexagon axioms for a multiplicative generating object $V_{\rm fund}$. }
Extending the coincidence of the two braided tensor structures from the special pairs or triple
of objects described in (c) of theorem \ref{Finkelberg_HL} to {\it all} irreducible objects   is    obstructed by   fusion multiplicities
and non-semisimplicity of $U_q({\mathfrak g})$. To overcome this and prove part (d), we implement a constructive cohomological reduction based on the abstract uniqueness results developed in    Theorems \ref{claim0} and Theorem \ref{claim1}. 

Theorem \ref{claim1} is a necessary   generalization of Theorem  \ref{claim0} to accomodate the bounds of our direct verification. 
 By evaluating the cohomology of the pentagon equation and its interplay with the hexagon equations, the theorem ensures that a braided tensor structure on a semisimple category is uniquely determined by its values on a generating object. This algebraic reduction eliminates the need for further direct evaluation on arbitrary modules.
 However, these theorems require a property that has not been considered so far: that the braid group representation between tensor powers of $V_{\rm fund}$ be generating. This links the problem  
of extending the equivalence to all objects to the field of generalized Schur-Weyl duality, ensuring that the braid group representations on tensor powers of $V_{\rm fund}$ uniquely determine the global braided tensor structure, thereby allowing us to identify the two categories based on their agreement on fusion of the generating object $V_{\rm fund}$.

 \section{The analytic Drinfeld twist versus the original Drinfeld--Kohno theorem}\label{10}\bigskip
  We next outline    general ideas   to prove part (b), (c) of Theorem \ref{Finkelberg_HL}. These parts concern the results on the 
  Drinfeld twist method for the Zhu algebra $A(V_{{\mathfrak g}_k})$ and the comparison of special braiding and associativity morphisms with those derived from Huang-Lepowsky theory. We comment below of analogies and differences with the original Drinfeld-Kohno theorem \cite{Drinfeld_quasi_hopf}.
  
  \subsection{Divergence: The use of the fundamental representation $V_{\rm fund}$, following \cite{Wenzl}}
  In  \cite{Wenzl} a fundamental generating representation $V_{\rm fund}$ is described for each simple Lie algebra ${\mathfrak g}$, and it is used
 for     the construction of the unitary structure in the associated quantum group fusion category.  
 In our approach we follow the work of Wenzl on the use of  the generating representation  $V_{\rm fund}$, both for ${\mathcal C}({\mathfrak g}, q)$ and ${\rm Rep}(A(V_{{\mathfrak g}_k})$. The common fusion structure of this representation allows to understand the structure
 of the associativity morphisms directly for the triples of objects described in part (c) of Theorem \ref{Finkelberg_HL} for
 the tensor structure of ${\rm Rep}(A_W({\mathfrak g}, q))$ and ${\rm Rep}(A(V_{{\mathfrak g}_k}))$ with respect to both tensor structures, arising from the twist construction or from Huang-Lepowsky tensor structure, by the results of Section \ref{8}. Having established the explicit unitary transport of the braided tensor structure
of $A_W({\mathfrak g}, q)$
 to the Zhu algebra via the fundamental representation $V_{\rm fund}$, we now compare our approach with the   work of Drinfeld. 
 This comparison highlights both  similarities with Drinfeld's work on coboundary quasi-Hopf algebras, and the novel aspects.
 
 \subsection{The intersection: The   twist as a square root of the coboundary} The deepest structural parallel between  Drinfeld's original work and our strategy lies in the explicit construction of a twist via a coboundary matrix. In his foundational paper on quasi-Hopf algebras \cite{Drinfeld_quasi_hopf}, Drinfeld introduces   \textit{coboundary} quasi-Hopf algebras: he provides a   formula for a twist by taking a \textbf{formal square root} of the coboundary matrix $\overline{R}$ of $U_h({\mathfrak g})$ over the ring of formal power series $\mathbb{C}[[h]]$. Our construction of the Drinfeld twist $T$ (Theorem 30.4 in \cite{CGP}) of $A_W({\mathfrak g}, q)$
 is the exact  operator-algebraic analogue of this formal square root. Because our framework operates at roots of unity and positive integer levels--where formal power series deformations are inapplicable--we cannot rely on formal Taylor expansions. Instead, we work with the analytic topology of the $C^*$-tensor category. By proving that the coboundary matrix $\bar{R}$ is strictly \textbf{self-adjoint} in this setting (Theorem 28.9 in \cite{CGP}), we rigorously define our twist $T$ as its analytic square root via functional calculus. This allows us to explicitly manage the braided symmetry on the special tensor products involving $V_{\rm fund}$, just as Drinfeld managed his formal category globally.
 
 The role of the associator in our work is a direct analogue of the Drinfeld-Kohno theorem. Just as Drinfeld uses a twist to relate a trivial associator to one derived from the KZ equations, we use the unitary twist $T$ to relate the trivial associator of the weak Hopf algebra $A_W(\mathfrak{g}, q)$ to the non-trivial Huang-Lepowsky associator (which is itself a VOA-theoretic analogue of the KZ structure).

\subsection{Divergence: Cohomological reduction and generalized Schur-Weyl duality} 
The divergence lies in the method of proof for the $2$-coboundary nature of the associator provided by the twist $T$.

 Global identity via Theorem  \ref{claim1}: Rather than a direct global   comparison, we use a tool of cohomological reduction. We verify the identity of the twisted   associator and the HL associator only on special pairs or triples of objects containing
 the generating representation $V_{\rm fund}$. Theorem  \ref{claim1} then provides   categorical methods via the pentagon and hexagon equations to uniqueley propagate this identity to the entire category.
 
 This theorem is based on validity of Schur-Weyl duality, a central notion in the theory of quantum groups that does not seem to appear
 in the original Drinfeld-Kohno theorem.

\subsection{The nature of the unitary structure} This analytic $C^*$-setting also reveals a striking property regarding the unitary structure itself. The unitary structure induced by the coboundary matrix exhibits profound differences from those of standard Hopf $C^*$-algebras (such as compact quantum groups, where the $^*$-involution commutes with the coproduct). In our setting of a unitary coboundary weak quasi-Hopf algebra, the $^*$-involution \textit{anticommutes} with the comultiplication. This property strongly resembles the antimultiplicativity of the involution in noncommutative topological spaces, and is intimately connected to preserving the balanced ribbon structure of the equivalence.

\subsection{Divergence: Braid group duality as a uniqueness tool}
Our approach diverges from the original Drinfeld--Kohno framework  of \cite{Drinfeld_quasi_hopf} by substituting Lie algebra cohomology with {\it braid group duality} as the fundamental mechanism for identifying the braided tensor structures.
This shift is necessitated by the rigidity obstruction at positive levels.

The core of our strategy relies on the property of \textit{multiplicative generation}: for the Lie types under consideration, the representation of the braid group in the centralizer algebras of the tensor powers of $V_{\rm fund}$ is ``large enough'' to uniquely determine the global associativity from the braiding data. As detailed in the following section, this allows us to use the de-quantization curve to establish a direct, algebraic path to the Finkelberg equivalence, bypassing the character-based verifications of earlier works. 

This conceptual framework provides the foundation for the technical inductive decomposition and the final identification proof presented in Sections \ref{11} and \ref{12}.

\section{Cohomological reduction of braided tensor equivalence to braid group duality for the proof of  Th. \ref{Finkelberg_HL}(\lowercase{d})}
\label{11}

To complete the proof of the ribbon braided tensor equivalence, we must transition from the analytic construction of the twist $T$ to a  result of structural identity of the two braided tensor structures, one arising from the Drinfeld twist method and the other from Huang-Lepowsky theory on the module category of the Zhu algebra. This section outlines the strategy used to propagate the identity of structures--initially verified on fundamental building blocks--to the entire category.

\subsection{The method of cohomological reduction}

The primary innovation in our identification strategy is the application of \textit{cohomological reduction}. Rather than attempting a global   comparison of the two categories, which meets significant obstructions at positive levels $k$ due to multiplicity decomposition, we utilize two abstract uniqueness results.  

These results (specifically Theorems  \ref{claim0} and  \ref{claim1}) establish that in a semisimple pre-tensor category, the global braided tensor structure is uniquely determined by its values on a specific set of \textit{special variables} provided a duality property between fusion tensor powers of the fundamental object $V_{\rm fund}$ and the braid group representation holds. By adopting this framework, we simplify the problem: if the structures coincide on a set of generating objects, the pentagon and hexagon equations force them to coincide globally.

\subsection{Identification on special pairs and triples}

We implement this reduction of the equality of associativity/braiding morphisms via the pentagon/hexagon by verifying the coincidence of our twisted braided tensor structure with the Huang--Lepowsky braided tensor  structure on the following ``building blocks'':
\begin{enumerate}
    \item \textbf{Special braiding morphisms:} The braided symmetries are verified on pairs where at least one variable is the fundamental representation $V_{\rm fund}$. Our Drinfeld-Kohno theorem was designed to this aim.
    \item \textbf{Special associativity morphisms:} The associativity is initially established on special triples where two   coordinates are fixed as $V_{\rm fund}$. The underlying strategy relies on reducing the equality of the associativity morphisms to purely extremal parenthesizations. (These special associativity morphisms are verified in our main application using  arguments analogous to the initial case.)

\end{enumerate}

The validity of this comparison in the application rests on the analytic properties of the \textit{unitary Drinfeld twist $T$}. As a non-commutative Radon--Nikodym derivative, $T$ ensures that the initial terms of the VOA primary fields are isometrically mapped to the corresponding structures in the quantum group category. For a wider explanation we refer the reader to Sects 8 and 37 of \cite{CGP}.

\subsection{The abstract uniqueness theorems}
\begin{defn} Let ${\mathcal C}$ be a (complex linear semisimple) category. Following Definition 4.1 in \cite{CGP}, we call ${\mathcal C}$ a {\it pre-tensor category} if it is equipped with a bilinear tensor product bifunctor $\otimes:{\mathcal C}\times{\mathcal C}\to{\mathcal C}$ and a tensor unit $\iota$ satisfying the usual unit axioms. We assume that ${\mathcal C}$ has subobjects and direct sums.

By Sect. 2.9 of \cite{EGNO} extended to pre-tensor categories, we may and shall assume that $\iota$ is a strict tensor unit: $\iota\otimes\rho=\rho\otimes\iota=\rho$ for every object $\rho$, and similarly for morphisms.

 We say that an object $V$ is {\it generating} if every irreducible of ${\mathcal C}$ is a subobject of some (parenthesized) tensor power of $V$.

A (complex, semisimple) {\it tensor category} is a pre-tensor category endowed with natural invertible associativity morphisms 
$\alpha_{\rho, \sigma, \tau}: (\rho\otimes\sigma)\otimes\tau\to\rho\otimes(\sigma\otimes\tau)$ satisfying the pentagon equation\begin{equation}\label{pentagon_equation}
\xymatrix@C=1em{
  ((\rho\otimes \sigma)\otimes \tau)\otimes \upsilon\ar[d]_{\alpha}\ar[r]^{\alpha\otimes {1}}
  &(\rho\otimes(\sigma\otimes \tau))\otimes \upsilon\ar[r]^{\alpha}
  &\rho\otimes((\sigma\otimes \tau)\otimes \upsilon)\ar[d]^{{1}\otimes\alpha}\\
  (\rho\otimes \sigma)\otimes(\tau\otimes \upsilon)\ar[rr]_{\alpha}&&\rho\otimes(\sigma\otimes (\tau\otimes \upsilon))
}
\end{equation} 
\end{defn}
We always assume $\alpha_{\rho, \iota, \tau}=1_{\rho\otimes\tau}.$

 \begin{defn} \label{braid_group_generating_property} Let 
  $({\mathcal C}, \otimes, \iota)$ be a semisimple linear pre-tensor category  
 and let $$c(\rho, \sigma): \rho\otimes\sigma\to\sigma\otimes\rho$$ be an invertible natural transformation
 such that $c(\iota, \sigma)$ and $c(\rho, \iota)$ are identity morphisms for all objects $\rho$ and $\sigma$.
 
 We shall say that an object $V$ of ${\mathcal C}$ satisfies {\it the   duality property with respect to $c$} if  
 for any positive integer $r\geq1$ the morphism space $(V^r, V^r)$ between any two   tensor powers of $V$   
   is linearly generated by finite compositions of  the component morphisms of $c$, their tensor products with identity morphisms and tensor products of morphisms in $(V, V)$ with identity morphisms.

 \end{defn}

\begin{defn}\label{braided_symmetry_tensor_category}
A {\it braided symmetry} for a tensor category $({\mathcal C}, \otimes, \iota, \alpha)$  is a natural isomorphism $c(\rho, \sigma)\in(\rho\otimes\sigma, \sigma\otimes\rho)$ satisfying
the normalization property
\begin{equation}\label{normalization_symmetry} c(\rho,\iota)=c(\iota,\rho)=1_\rho\end{equation} and such that the following two {\it hexagonal diagrams} commute
\begin{equation}\label{braided_symmetry1} \begin{CD}
 (\rho\otimes \sigma)\otimes \tau @>{{\alpha}}>>  \rho\otimes(\sigma\otimes \tau)@>{c}>> (\sigma\otimes\tau)\otimes \rho \\
@V{c\otimes 1}VV @. @VV{\alpha}V \\
 (\sigma\otimes\rho)\otimes\tau@>{\alpha}>>  \sigma\otimes(\rho\otimes\tau)@>{1\otimes c}>> \sigma\otimes(\tau\otimes\rho) \end{CD}\end{equation}
 \begin{equation}\label{braided_symmetry2} \begin{CD}
 (\rho\otimes \sigma)\otimes \tau @>{{c}}>>  \tau\otimes(\rho\otimes \sigma)@>{\alpha^{-1}}>> (\tau\otimes\rho)\otimes\sigma \\
@A{\alpha^{-1}}AA @. @AA{c\otimes 1}A \\
 \rho\otimes(\sigma\otimes\tau)@>{1\otimes c}>>  \rho\otimes(\tau\otimes\sigma)@>{\alpha^{-1}}>> (\rho\otimes\tau)\otimes\sigma \end{CD}\end{equation}
 \end{defn}

    \begin{thm}\label{claim0}
 Let $({\mathcal C}, \otimes, \iota)$ be a semisimple pre-tensor category with  a generating object $V$ and admitting a faithful weak quasi-tensor functor $({\mathcal F}, F, G): {\mathcal C}\to{\rm Vec}$ into the category of finite dimensional vector spaces.
 
  Let  $c(\rho, \sigma): \rho\otimes\sigma\to\sigma\otimes\rho$ be an invertible natural transformation satisfying
  $c(\iota, \rho)=c(\rho, \iota)=1_\rho$ for all objects $\rho$.
  Let
  $V$ satisfy the   duality property with respect to $c$.

 Let $\alpha$ and $\beta$ be two associativity morphisms for $({\mathcal C}, \otimes, \iota)$ such that $({\mathcal C}, \otimes, \iota, \alpha, c)$ and $({\mathcal C}, \otimes, \iota, \beta, c)$ are braided tensor categories.
 Set
\begin{equation}\label{definition_of_V}
 {\mathcal V}=\{(V_\lambda, V, V), (V, V_\lambda, V), (V, V, V_\lambda), \quad V_\lambda\in{\rm Irr}({\mathcal C})\}.\end{equation}

Assume that     
 \begin{equation}\label{assumption_on_associators2}
 \alpha=\beta \quad\text{on }
 {\mathcal V}.\end{equation}  
  Then  $\alpha=\beta$ everywhere.  
 \end{thm}

  \begin{proof}

a)    By naturality of the associativity morphisms, if $\alpha_{V_1, V_2, V_3}=\beta_{V_1, V_2, V_3}$
  for a fixed triple $(V_1, V_2, V_3)$ then $\alpha_{W_1, W_2, W_3}=\beta_{W_1, W_2, W_3}$ for any other triple $(W_1, W_2, W_3)$ such that each $W_i$ is equivalent to $V_i$.

 By   semisimplicity of the category,   the fact that $V$ is generating, and naturality of the associators in the three variables, it follows that the assumptions in the statement regarding the complete fixed collection   $\{ V_\lambda\}$ of irreducible objects
is equivalent to the corresponding assumptions   regarding a collection of one (and hence all) fixed choice of paranthesizations in a tensor power
 $V^r$ of $V$ for each order $r\geq 1$. Similarly, the thesis of the statement is equivalent to
showing that $\alpha_{ V^s, V^t, V^u  }=\beta_{ V^s, V^t, V^u  }$ with  $(V^s, V^t, V^u)$ a fixed
choice of   tensor powers of $V$   for each triple   $(s, t, u)$ of corresponding orders. 

We   then pass to this equivalent formulation of   the statement replacing the irreducible arbitrary term
$V_\lambda$ with an arbitrary tensor power $V^r$ for   each non negative order $r$. When necessary, we shall
make a choice of tensor powers $V^r$ of $V$  that we specify in the course of the proof.

To show that $\alpha_{V^s, V^t, V^u  }=\beta_{ V^s, V^t, V^u  }$ on a given triple of tensor powers, we may assume     $\min\{s, t, u\}\geq1$, as when one of $s$, $t$, or $u$ is zero, then
the two associativity morphisms  equal   the identity map, by definition. 
\medskip

b) On any triple $(s, t, u)$ of positive integers, we   define the integer valued function
$$f(s, t, u) := s+t+u.$$
\medskip     
c)    Consider also  the   integer-valued function $$g(s, t, u):=\min\{s, t, u\}+\text{medium}\{s, t, u\},$$  the sum of the minimum and medium value
    among the three variables $(s, t, u)$. 
    \medskip

d) We proceed by induction. First inductive assumption on $f$: 
\medskip

We have $f(s,t,u)\geq3$, and $f(s,t,u) =3$ if and only if  $(s, t, u)=(1,1,1)$. In this case,   $\alpha$ and $\beta$ coincide  on the corresponding triple  $(V, V, V)$ by assumption.
Let $N$ be an integer $\geq3$. Assume  that $\alpha_{V^s,V^t,V^u} = \beta_{V^s,V^t,V^u}$ for all triples $(s, t, u)$ such that $f(s, t, u) \leq N$. \medskip

e) First inductive step on $f$:

Let    
 $(V^s,V^t,V^u)$ be  a new triple with $$f(s,t,u) = N+1.$$
   Then at least one among $s$, $t$, $u$ takes a value $>1$. 
   \medskip

f) Second inductive assumption on $g$, with the value of $f$ fixed to $N+1$:
\medskip

The function $g$ takes minimum value $=2$, and this happens if and only if  two coordinates equal $1$ and the third is an arbitrary tensor power of $V$. On such triples, $\alpha$ and $\beta$ coincide by assumption.
Let $M$ be an integer $\geq2$.\medskip

 Let us consider
  triples $(s, t, u)$ on which    $f$  takes the value $N+1$.
     Assume   that for each such $(s, t, u)$, we   have
    $$\alpha_{V^s,V^t,V^u} = \beta_{V^s,V^t,V^u}, \quad \text{for }
    g(s, t, u)\leq M.$$\medskip

 \medskip

g) Second inductive step on $g$ (at $f=N+1$): 
Let $(s, t, u)$ be a new  triple for which $f(s, t, u)=N+1$ and $g(s, t, u)=M+1$.
Then at least two coordinates among $(s, t, u)$ take a value $>1$.

\medskip

h) We  
need to show   that $\alpha=\beta$ in the second inductive step for $g=M+1$ (at $f=N+1$).

We do this in the following paragraphs from i) to p).
 \medskip

i) The maximum coordinate in $(s, t, u)$ in the inductive step on $g$ takes the value $N+1-(M+1)=N-M$.
    
If   $(s, t, u)$ satisfies the second inductive assumption
    on $g$, then the maximum coordinate in $(s, t, u)$ takes a value $\geq N-M+1$, strictly larger than the maximum coordinate    in the inductive step on $g$ (at $f$ fixed to the value $N+1$). This property distinguishes the second inductive assumption on $g$ (at $f$ fixed to the value $N+1$.)

We next analyse
 case by case what happens to the associativity morphisms $\alpha_{V^s, V^t, V^u}$ on each coordinate
 $s$, $t$, or $u$    $>1$, using the pentagon equation.
Notice that the following equations are symmetric in  $s$ and $u$, but not in $t$.

 Case $t>1$. In this case, we write
   $$V^t=V^{t_1}\otimes V^{t_2}$$ according to the parenthesization of $V^t$. We have
   $t_1$, $t_2\geq1$. In the following we use the notation
   $$V^t=V^{t_1+t_2}.$$
   Then the pentagon equation  (\ref{pentagon_equation}),    gives
 \begin{equation}\label{dec1} 
 \alpha_{V^s, V^t, V^u}=1_s\otimes\alpha_{V^{t_1}, V^{t_2}, V^u}^{-1}\circ \alpha_{V^s, V^{t_1}, V^{t_2+u}}\circ \alpha_{V^{s+t_1}, V^{t_2}, V^u}\circ \alpha_{V^s, V^{t_1}, V^{t_2}}^{-1}\otimes 1_u.
   \end{equation}
   
     Case $s>1$. 
   We write  $V^s=V^{s_1}\otimes V^{s_2},$ $s_1$, $s_2\geq1$, according to the given paranthesization of $V^s$.
 Then by the pentagon equation  (\ref{pentagon_equation}),
 \begin{equation}\label{dec2} \alpha_{V^s, V^t, V^u}=\alpha^{-1}_{V^{s_1}, V^{s_2}, V^{t+u}}\circ 1_{s_1}\otimes \alpha_{V^{s_2}, V^{t}, V^u}\circ\alpha_{V^{s_1}, V^{s_2+t}, V^u}\circ\alpha_{V^{s_1}, V^{s_2}, V^t}\otimes 1_u.
    \end{equation}
    
      Case $u>1$.  We  write  
$V^u=V^{u_1}\otimes V^{u_2}$ with $u_1$, $u_2\geq1$.
By the pentagon equation (\ref{pentagon_equation}),
 \begin{equation}\label{dec3} \alpha_{V^s, V^t, V^u}=1_s\otimes \alpha_{V^t, V^ { u_1 }, V^ { u_2 } } \circ\alpha_{V^s,  V^{ t+u_1}, V^ { u_2 }}\circ\alpha_{V^s, V^t, V^ { u_1}}\otimes 1_{u_2 }\circ\alpha^{-1}_{V^{s+t   },  V^ {u_1}, V^ {u_2}}.
     \end{equation}
     
     We next assume the hypotheses in second inductive step on $g$. As remarked, a 
       second coordinate in $(s, t, u)$ is automatically $>1$.
     
     We divide the proof of the second inductive step on $g$ in the following sub-cases $A_1$ and $A_2$.

     Sub-case $A_1$:
$t$ is maximum, so $t>1$. 
      We assume that the second coordinate $>1$    is $s$ and we use (\ref{dec2}).
     (If instead $u>1$, we may  reason in a similar way using (\ref{dec3}),  by the symmetry between (\ref{dec2}) and (\ref{dec3}).)

    Then  $\alpha$ and $\beta$ coincide on the second and the fourth factor of the right hand side of (\ref{dec2}) by the first inductive assumption on $f$. They also coincide on the first and third factor of the right hand side of (\ref{dec2}) by the second inductive assumption on $g$, since
      the corresponding triples of those factors satisfy the first inductive step  ($f$  still takes the value $N+1$) and the maximum of those triples  is $t+u$ and $s_2+t$ respectively, which it is $> t$ in both cases. 
      
      We conclude that $\alpha=\beta$ on triples $(s, t, u)$ satisfying the assumptions of the second inductive step on $g$, with the restrictive assumption that $t$ is   maximum.
      The proof of   sub-case $A_1$   is complete.

     Sub-case $A_2$:
   the maximum is $s$ or $u$ but not $t$. This case in turn divides in two symmetric sub-cases:
     
  $(A_2)_1$   $u$ is maximum, but $t$ is not maximum, so $u>1$ and $t<u$. We also have $s\leq u$.
     
    $(A_2)_2$  $s$ is maximum but $t$ is not maximum, so $s>1$, $t<s$. We also have $u\leq s$.
     
It suffices to  show   $(A_2)_1$ by the symmetry of the pentagon equations (\ref{dec2}) and (\ref{dec3}) on the first and last coordinate.

The case $(A_2)_1$ in turn divides in the following sub-cases $((A_2)_1)_1$ and $((A_2)_1)_2$.

Subcase $((A_2)_1)_1$ of   $(A_2)_1$:   $s=u$. Then $s>1$. We may decompose $V^s$ using (\ref{dec2}) again. We see that $\alpha$ and $\beta$ coincide
       again on the second and fourth factor at the right hand side by the inductive assumption on $f$ and on the first factor by the inductive assumption
       on $g$, since on this factor the maximum is $t+u>u$. Thus $\alpha$ and $\beta$ coincide on $(s, t, u)$ with $u$ maximum, $t<u$, $s=u$,
       $f(s, t, u)=N+1$ and $g(s, t, u)=M+1$ if and only if $\alpha$ and $\beta$ coincide on the triple $(s_1, s_2+t, u)$. We study coincidence
       of the associativity morphisms on the last triple.
       We may assume $s_2+t=u$ by naturality of the associators. Then on $(s_1, s_2+t, u)$ $f$ takes the same value as on $(s, t, u)$, 
       the maximum is $s_2+t= s=u$ and    $g$ takes the value
       $s+t=M+1$. Thus we are in the second inductive step 
       with maximum in the middle. This case has been completed above, and says that  $\alpha=\beta$ on $(s_1, s_2+t, u)$ as well, and it follows that the proof of   sub-case $((A_2)_1)_1$ of   $(A_2)_1$ is complete.

Sub-case $((A_2)_1)_2$ of $(A_2)_1$:     $s<u$. Recall that   the maximum in $(s, t, u)$ is $u>1$, $t<u$.  
If $s=1$ we shall not proceed for the moment. If $s>1$ then we proceed with the aim to  reduce to the case $s=1$ as follows.

   We apply (\ref{dec2}).      By the same arguments as in $((A_2)_1)_1$, if for some $s_2$ such that $s=s_1+s_2$ we have $s_2+t\geq u$ then we obtain equality of $\alpha$ and $\beta$ on this triple.
   
   We are left to analyse equality of the associativity morphisms
 in the case $s_2+t< u$ for all possible decompositions of $s=s_1+s_2$.
      Thus we choose the extremal decomposition with $s_1=1$ and $s_2=s-1$. We set $\ell=s_2+t=s+t-1$.
       
In conclusion analysis on all possible cases reduces
 to   triples of the form $$(1, \ell, u), \quad\quad \ell<u.$$ (case $(A_2)_1$) and of the form
       $$(s, m, 1), \quad\quad m<s$$
       (case $(A_2)_2$.)

      By the assumption in   the second inductive step,
         $\ell=m=M<u$ and $u=s=N-M$. Thus we are reduced to study equality of the associativity morphisms $\alpha$ and $\beta$ on
         the objects determined determined by the two inductive steps corresponding to the triples
         $$(1, \ell, u), \quad (u, \ell, 1), \quad \quad \ell=M, \quad u=N-M.$$
         Passing to the inverse of the associativity morphisms on the second triple, we need to show that 

  $$\alpha=\beta \quad\text{on}\quad (V\otimes V^\ell)\otimes V^u, \quad \alpha^{-1}=\beta^{-1} \quad \text{on}\quad V^u\otimes (V^\ell\otimes V),$$

 Equation (\ref{dec1}) with $t_2=1$ applied to $(1, \ell, u)$    says that studying equality of  the associativity morphisms $\alpha$ and $\beta$ on $(1, \ell, u)$ is equivalent to
       studying the same question on $(\ell, 1, u)$. Similarly the case $(s, m, 1)$ with $m<s$ is equivalent to $(s, 1, m)$.\medskip
       
 l)      So we   consider the following four  objects of ${\mathcal C}$ of  order $N+1$,  with $\ell=M<u=N-M$, so $M<N/2$,
       
    \begin{equation}\label{reduction} (V\otimes V^\ell)\otimes V^u, \quad (V^\ell\otimes V)\otimes V^u, \quad V^u\otimes (V^\ell\otimes V), \quad V^u\otimes (V\otimes V^\ell).    \end{equation}

        For convenience of the proof, we include each of the four objects in (\ref{reduction}) in a  larger finite family that also includes  triples for which equality of the two associativity morphisms is known, by the second inductive assumption on $g$.

In order from left to right, each of the above four objects belongs to the following corresponding      family of objects of ${\mathcal C}$, from 1) to 4), respectively:  \medskip
  
    1)   $(V\otimes V^\ell)\otimes V^u$, with $\ell\leq M$, $\ell<u$, $\ell+u= N$,
  
  2) $(V^\ell\otimes V)\otimes V^u$, with $\ell\leq M$, $\ell<u$, $\ell+u= N$,

    3)  $V^u\otimes (V^\ell\otimes V)$, with $\ell\leq M$, $\ell<u$, $\ell+u= N$,
  
  4)   $V^u\otimes (V\otimes V^\ell)$, with $\ell\leq M$, $\ell<u$, $\ell+u= N$.

    The family   1) intersects  2)   if and   only if $\ell=1$ (we already know that $\alpha=\beta$
    by assumption in this case, but we shall not use this in the following, but only at the end).
    
    The family 1) intersects 3) if and only il only if $\ell=1$   and $u=2$.

    The family 1) intersects 4) if and only if $\ell+1=u$ and $V^u=V\otimes V^\ell$.
    
     The family 2) intersects 3) if and only if $\ell+1=u$ and $V^u=V^\ell \otimes V$.

    The family 2) intersects 4) if only if $\ell=1$    and $u=2$.
         
 The family 3) intersects 4)  if and only if $\ell=1$ and $u=2$.

  We define the following morphisms in ${\mathcal C}$ between objects of 1), 2), 3), 4), respectively:
 
 $$ T_{(V\otimes V^\ell)\otimes V^u}:=(\beta^{-1}\alpha)_{V, V^\ell,  V^u}: (V\otimes V^\ell)\otimes V^u\to (V\otimes V^\ell)\otimes V^u,$$
 
 $$T_{(V^\ell\otimes V)\otimes V^u}:=(\beta^{-1}\alpha)_{V^\ell, V,  V^u}: (V^\ell\otimes V)\otimes V^u\to (V^\ell\otimes V)\otimes V^u,$$

       $$ T_{V^u\otimes (V^\ell\otimes V)}:=(\beta\alpha^{-1})_{V^u, V^\ell,  V}: V^u\otimes (V^\ell\otimes V)\to V^u\otimes (V^\ell\otimes V).$$
       
       $$ T_{V^u\otimes (V\otimes V^\ell)}:=(\beta\alpha^{-1})_{V^u, V,  V^\ell}: V^u\otimes (V\otimes V^\ell)\to V^u\otimes (V\otimes V^\ell),$$

We next see that on objects belonging to the intersections, these morphisms are defined in the same way. This is clear  for the intersections 
$1)\cap 2)$,  $3)\cap 4)$ by definition of the $T$-morphisms. 

We next consider     the intersections   $1)\cap 4)$ and $2)\cap 3)$.
The intersection $1)\cap 4)$ is the object $(V\otimes V^\ell)\otimes (V\otimes V^\ell)\in 1)\cap 4)$.
We need to show that
$$(\beta^{-1}\alpha)_{V, V^\ell,  V\otimes V^\ell}=(\beta\alpha^{-1})_{V\otimes V^\ell, V,  V^\ell}$$
equivalently that
$$\beta_{V,  V^\ell, V\otimes V^\ell}\circ \beta_{V\otimes V^\ell, V, V^\ell}=   
\alpha_{V, V^\ell,  V\otimes V^\ell}\circ\alpha_{V\otimes V^\ell, V,  V^\ell}.$$
We use the pentagon equation (\ref{pentagon_equation}) for $\alpha$ and $\beta$ and   the above equation  
becomes
$$ 1\otimes\beta\circ  \beta_{V, (V^\ell\otimes V),  V^\ell}\circ\beta\otimes 1= 1\otimes\alpha\circ  \alpha_{V, (V^\ell\otimes V),  V^\ell}\circ\alpha\otimes 1.                         $$
The first and third associativity morphisms in $\beta$ from on the left hand side of the equation  equal the corresponding associativity morphisms in $\alpha$ at the right hand side, by the first inductive assumption on $f$.
The middle associativity morphisms on both sides  are also equal   as we have already completed the proof of the second inductive step on $g$ in the case where
 the middle term has maximal power.
 Thus  the above $T$-morphisms coincide on $1)\cap 4)$. One similarly shows that they coincide on $2)\cap 3)$.
 
  The intersections $1)\cap 3)$, $2)\cap 4)$ are a particular case of the last two intersections for $\ell=1$.
  This example is indeed quite instructive.
 
  Let us consider the linear full subcategory ${\mathcal D}$ of   ${\mathcal C}$ with objects $1)\cup 2)\cup 3)\cup 4)$.
  This subcategory depends on $N$ and $M$, and therefore the same holds for further constructions associated to
  ${\mathcal D}$ in the rest of the proof.

   We next see that the $T$-morphisms are the components of  a natural transformation $T$ from the identity functor of ${\mathcal D}$
  to itself. From this property it will follow that $T$ extends uniquely to  the completion $\tilde{\mathcal D}$ of ${\mathcal D}$ with subobjects  and direct sums.
  
  We need to show    the natural transformation property of $T$, i.e. that if $X$ and $Y$ are objects in ${\mathcal D}$
  and $S: X\to Y$ is a morphism in ${\mathcal C}$ then
  $$S\circ T_X=T_Y\circ S.$$
  To show this property, 
we need to know the morphisms $S$ in ${\mathcal C}$, or at least a generating family under composition and linear combination.
These morphisms are described by the braid group generating property of the braided symmetry $c$ in ${\mathcal C}$ by assumption.
Recall that we do not need to specify associativity morphisms to describe this property.
\medskip

m)  In defining ${\mathcal D}$, we let the factors $V^\ell$ and $V^u$ appearing  in the tensor product  objects as in 1), 2), 3), 4), vary in the same orbit under the action of the braid group.
  By naturality of the associativity morphisms $\alpha$ and $\beta$ in the variables, this gives in case 1)
  commutation relations of the kind
  $$T\circ (c_1\otimes c_2)\otimes c_3=(c_1\otimes c_2)\otimes c_3\circ T$$
  where $c_1: V\to V$, $c_2: V^\ell\to V^\ell$, and $c_3: V^u\to  V^u$ are  morphisms in $(V, V)$ or braiding morphisms or their tensor products    with identity morphisms.\medskip

n) We next use the hexagonal diagrams (\ref{braided_symmetry1})  and (\ref{braided_symmetry2}) for $(\alpha, c)$ and $(\beta, c)$
to complete   the   generators of the morphisms in  ${\mathcal D}$ and the corresponding
  commutation relations. This will   complete the proof that $T$ is a natural transformation.

These hexagonal diagrams   allow us
  to put the maximal power  $V^u$ of the triple in the middle. In this case we have already shown the second inductive step on $g$, so we know that
    two associators $\alpha$ and $\beta$ (or their inverses) coincide. We aim
   to extract further information  from this, for the   cases cases where $V^u$ is not in the middle.

In the hexagonal diagram  (\ref{braided_symmetry1}) we start choosing
  $$\rho=V, \quad \tau=V^u, \quad \sigma=V^\ell.$$
 Since   the object $\tau=V^u$   with maximum power is in the middle place on the right vertical arrow of the diagram,
 we   know
 that the corresponding associativity morphisms coincide $\alpha=\beta$ coincide  on the
 triple $(V^\ell, V^u, V)$.   
 Thus after equating  these known coinciding  associativity morphisms and the braiding morphism at their left and right in the corresponding equation in  (\ref{braided_symmetry1}),
 we obtain   the following equality of the other hand sides,
  $$\alpha_{V^\ell, V,  V^u}\circ c(V, V^\ell)\otimes 1\circ \alpha^{-1}_{V, V^\ell, V^u} =\beta_{V^\ell, V,  V^u}\circ c(V, V^\ell)\otimes 1\circ \beta^{-1}_{V, V^\ell, V^u},$$
 so
 $$(\beta^{-1}\alpha)_{V^\ell, V,  V^u}\circ c(V, V^\ell)\otimes 1 =  c(V, V^\ell)\otimes 1\circ (\beta^{-1}\alpha)_{V, V^\ell, V^u},$$
  therefore
    \begin{equation}\label{first_commutation_with_braiding_T} T_{(V^\ell\otimes V)\otimes V^u}\circ c(V, V^\ell)\otimes 1 =  c(V, V^\ell)\otimes 1\circ T_{(V\otimes V^\ell)\otimes V^u}.
    \end{equation}

  After writing  (\ref{braided_symmetry2}) in the same form as  (\ref{braided_symmetry1}) for the inverse braided symmetry
  $c'(\rho, \sigma):=c(\sigma, \rho)^{-1}$, diagram (\ref{braided_symmetry2}) similarly gives

  $$(\beta^{-1}\alpha)_{V, V^\ell, V^u}\circ c(V^\ell, V)\otimes 1=c(V^\ell, V)\otimes 1\circ(\beta^{-1}\alpha)_{V^\ell, V,  V^u},$$
   that is
      \begin{equation}\label{second_commutation_with_braiding_T} T_{(V\otimes V^\ell)\otimes V^u}\circ c(V^\ell, V)\otimes 1=c(V^\ell, V)\otimes 1\circ T_{(V^\ell\otimes V)\otimes V^u}.
    \end{equation}

  Next, we use the top left associativity morphism in (\ref{braided_symmetry1}) with 
 $$\sigma=V^u, \quad \rho=V, \quad \tau=V^\ell.$$
 With this choice the two top left associativity morphisms $\alpha$ and $\beta$
 in the corresponding hexagonal diagrams coincide, because the maximal power $V^u$    is in the middle.
   
  Equating the  hand sides of these known coinciding braiding and associativity morphisms similarly gives the following equation of the other hand sides,
 $$\alpha^{-1}_{V^u, V^\ell, V}\circ1\otimes c(V, V^\ell)\circ \alpha_{V^u, V, V^\ell}=\beta^{-1}_{V^u, V^\ell, V}\circ1\otimes c(V, V^\ell)\circ \beta_{V^u, V, V^\ell},$$
 so
  $$(\beta\alpha^{-1})_{V^u, V^\ell, V}\circ1\otimes c(V, V^\ell) = 1\otimes c(V, V^\ell)\circ (\beta\alpha^{-1})_{V^u, V, V^\ell},$$   that is
     \begin{equation}\label{third_commutation_with_braiding_T} T_{V^u\otimes(V^\ell\otimes V)}\circ1\otimes c(V, V^\ell) = 1\otimes c(V, V^\ell)\circ T_{V^u\otimes(V\otimes V^\ell)}.
    \end{equation}
 
 Similarly (\ref{braided_symmetry2}) gives
$$(\beta\alpha^{-1})_{V^u, V, V^\ell}\circ1\otimes c(V^\ell, V)= 1\otimes c(V^\ell, V)\circ(\beta\alpha^{-1})_{V^u, V^\ell, V},$$
    so
    \begin{equation}\label{fourth_commutation_with_braiding_T}T_{V^u\otimes(V\otimes V^\ell)}\circ1\otimes c(V^\ell, V)= 1\otimes c(V^\ell, V)\circ T_{V^u\otimes(V^\ell\otimes V)}.
    \end{equation}
\medskip

 p) Finally we use the bottom left associativity morphism in (\ref{braided_symmetry1}) with
 $$\rho=V^u, \quad \sigma=V^\ell, \quad \tau= V.$$ The two bottom left associativity morphisms $\alpha$ and $\beta$ coincide again for the same reasons as above. 
 Putting the unknown associativity morphisms on the other side of this equation for the two cases $\alpha$ and $\beta$ and    equating the other hand sides gives
 
 $$\alpha_{V^\ell, V, V^u}\circ c(V^u, V^{\ell}\otimes V)\circ\alpha_{V^u, V^\ell, V}=\beta_{V^\ell, V, V^u}\circ c(V^u, V^{\ell}\otimes V)\circ\beta_{V^u, V^\ell, V},$$
 so
  $$(\beta^{-1}\alpha)_{V^\ell, V, V^u}\circ c(V^u, V^{\ell}\otimes V)=c(V^u, V^{\ell}\otimes V)\circ(\beta\alpha^{-1})_{V^u, V^\ell, V},$$
thus
  \begin{equation}\label{fifth_commutation_with_braiding_T} T_{(V^\ell\otimes V)\otimes V^u}\circ c(V^u, V^{\ell}\otimes V)=c(V^u, V^{\ell}\otimes V)\circ T_{V^u\otimes(V^\ell\otimes V)}.
  \end{equation}

Similarly (\ref{braided_symmetry2}) gives

 $$(\beta\alpha^{-1})_{V^u, V^\ell, V}\circ c( V^{\ell}\otimes V, V^u)=c(V^{\ell}\otimes V, V^u)\circ(\beta^{-1}\alpha)_{V^\ell, V, V^u},$$
  therefore
    \begin{equation}\label{sixth_commutation_with_braiding_T} T_{V^u\otimes(V^\ell\otimes V)}\circ c( V^{\ell}\otimes V, V^u)=c(V^{\ell}\otimes V, V^u)\circ T_{(V^\ell\otimes V)\otimes V^u}.
   \end{equation}
\medskip
 
 q)   The six commutation relations  (\ref{first_commutation_with_braiding_T}),  (\ref{second_commutation_with_braiding_T}) (\ref{third_commutation_with_braiding_T}) (\ref{fourth_commutation_with_braiding_T}) (\ref{fifth_commutation_with_braiding_T}) (\ref{sixth_commutation_with_braiding_T}) complete the  proof that $T$ is a natural transformation from the identity functor of
    ${\mathcal D}$ to itself.
    Thus $T$ extends as anticipated to a natural transformation on the completion  $\tilde{\mathcal D}$
    of  ${\mathcal D}$ with subobjects and direct sums.
    
    Let $\tilde{A}$ be the discrete finite dimensional algebra corresponding to the restriction of the linear functor ${\mathcal F}$ to
     $\tilde{\mathcal D}$. By linear Tannakian duality, $T$ corresponds to an invertible element of $\tilde{A}$ (which must be central).
     Since $A$ is discrete, there is a unique central partially invertible  element $U\in A$ with support in the semisimple quotient $\tilde{A}$
     such that
     $${\mathcal F}(T_X)={\mathcal F}(X)(U), \quad\quad X\in \tilde{\mathcal D}.$$
     Thus  for the objects in 1),
      $${\mathcal F}(T_{(V\otimes V^\ell)\otimes V^u})=\pi\otimes \pi^{\otimes \ell}\otimes \pi^{\otimes u}
     (\Delta\otimes 1\circ \Delta(U))\quad \ell\leq M, \quad \ell<u, \quad \ell+u=N,$$
     with $\Delta$ the coproduct of ${A}$ defined by Tannakian duality applied to $({\mathcal F}, F, G)$,
     $\pi$ the representation of $A$ corresponding to $V$ via duality, and $\pi^{\otimes k}$ a suitable tensor power
     of $\pi$ of order $k$ as a representation of the weak quasi-bialgebra $A$, corresponding to a given tensor power
     $V^k$.

Thus
    \begin{equation}\label{the_central_element}
         \beta_{V, V^\ell, V^u}=\alpha_{V, V^\ell, V^u}\circ \pi\otimes \pi^{\otimes \ell}\otimes \pi^{\otimes u}
     (\Delta\otimes 1\circ \Delta(U^{-1})),
      \end{equation}
      
      for $ \ell\leq M, \quad \ell<u, \quad \ell+u=N$.
      
      Similarly for the cases 2), 3), 4):
       \begin{equation}   \beta_{V^\ell, V, V^u}=\alpha_{V^\ell, V, V^u}  \circ \pi^{\otimes \ell}\otimes \pi\otimes \pi^{\otimes u}
     (\Delta\otimes 1\circ \Delta(U^{-1})),
         \end{equation}
         
         for $ \ell\leq M, \quad \ell<u, \quad \ell+u=N$,
         
        \begin{equation}      \beta_{V^u, V^\ell, V}= \pi^u\otimes \pi^{\otimes \ell}\otimes \pi
     (1\otimes \Delta\circ \Delta(U)) \circ  \alpha_{V^u, V^\ell, V},  
              \end{equation}
              
               for $ \ell\leq M, \quad \ell<u, \quad \ell+u=N$,
               
               and finally

                \begin{equation}      \beta_{V^u, V, V^\ell}= \pi^u\otimes\pi\otimes \pi^{\otimes \ell} 
                     (1\otimes \Delta\circ \Delta(U))  \circ \alpha_{V^u, V^\ell, V},
              \end{equation}
     
     for $ \ell\leq M, \quad \ell<u, \quad \ell+u=N$.
     
     It suffices to show that only in one of these cases, $U=1$ to obtain the claim of the inductive step on $g$, because $U$ is a common element for all four cases.
     
   The subset of objects of ${\mathcal D}$ in each of the families from 1) to 4) with $\ell=1$, defines a full linear subcategory with subobjects and direct sums  $\tilde{\mathcal E}_i$. Obviously $\tilde{\mathcal E}_1=\tilde{\mathcal E}_2$ 
  and  $\tilde{\mathcal E}_3=\tilde{\mathcal E}_4$. Both $\tilde{\mathcal E}_1$ and $\tilde{\mathcal E}_3$ are equivalent to $\tilde{\mathcal D}$. 
   
   Thus  $\tilde{\mathcal E}_1$ or $\tilde{\mathcal E}_3$ define each the same  quotient algebra  $\tilde{A}$.
   Since the coproduct of $A$ is faithful and since the two associators $\alpha$ and $\beta$ coincide
   on   objects of $\tilde{\mathcal E}_1$ or $\tilde{\mathcal E}_3$ by our assumption, it follows that  $U$ is the identity of ${\tilde{A}}$, and therefore $\alpha$ and $\beta$ coincide on all the objects of  $\tilde{\mathcal D}$.
   
       By definition this gives
     $$\alpha=\beta$$ on triples of the form
     $$(V\otimes. V^\ell)\otimes V^u, \quad (V^\ell\otimes V)\otimes V^u, \quad (V^u\otimes V^\ell)\otimes V, \quad (V^u\otimes V)\otimes V^\ell,$$
     for all $\ell\leq M$, $\ell<u$, $\ell+u=N$. 
     This completes the proof of the second inductive step for $g=M+1$ at $f=N+1$, stated in paragraph h).

  \end{proof}  
 
 \begin{defn}

Let us consider associativity morphisms for $n\geq 4$:
$$\alpha'_{W_1, W_2, W_3, \dots, W_n}: (\dots(W_1\otimes W_2)\otimes W_3)\otimes\dots )\otimes W_n \to W_1\otimes (\dots \otimes(W_{n-2}\otimes(W_{n-1}\otimes W_n)\dots)$$
that pass from left-parenthesized   to the right-parenthesized tensor products of the   objects $W_i$.
Let us also consider their inverses
$$(\alpha'_{W_1, W_2, W_3, \dots, W_n})^{-1}: 
W_1\otimes (\dots \otimes(W_{n-2}\otimes(W_{n-1}\otimes W_n)\dots)\to (\dots(W_1\otimes W_2)\otimes W_3)\otimes\dots )\otimes W_n$$
that pass from right-parenthesized   to the left-parenthesized tensor products of the   objects $W_i$.
We refer to such associativity morphisms  $\alpha'_{W_1, W_2, W_3, \dots, W_n}$   on $n$ variables and their inverses
$(\alpha'_{W_1, W_2, W_3, \dots, W_n})^{-1}$
 as {\it extremal}.
 \end{defn}

    \begin{thm}\label{claim1}
 Let $({\mathcal C}, \otimes, \iota)$ be a semisimple pre-tensor category with  a generating object $V$ and admitting a faithful weak quasi-tensor functor $({\mathcal F}, F, G): {\mathcal C}\to{\rm Vec}$ into the category of finite dimensional vector spaces.
 
 Let  $c(\rho, \sigma): \rho\otimes\sigma\to\sigma\otimes\rho$ be a normalized invertible natural transformation and let
  $V$ satisfy the   duality property with respect to $c$.
  
   Let  $d(\rho, \sigma): \rho\otimes\sigma\to\sigma\otimes\rho$ be another normalized invertible natural transformation.
   
    Let $\alpha$ and $\beta$ be two associativity morphisms for $({\mathcal C}, \otimes, \iota)$ such that $({\mathcal C}, \otimes, \iota, \alpha, c)$ and $({\mathcal C}, \otimes, \iota, \beta, d)$ are braided tensor categories.

Let ${\mathcal V}$ be defined as in (\ref{definition_of_V}).

Assume that for all 
$\lambda \in{\rm Irr}({\mathcal C})$: 

\begin{itemize} 

\item[(a)] $c(V_\lambda, V)=d(V_\lambda, V)$ and $c(V, V_\lambda)=d(V, V_\lambda)$ .

\item[(b)]
 \begin{equation}  
 \alpha=\beta \quad\text{on }
 {\mathcal V}.\end{equation}

 \item[(c)] Extremal associativity morphisms  $\alpha'$ and $\beta'$  coincide on   $1+k$-tuples of the form
 $(V^r, V, V, \dots, V)$ with $V$ repeated a number of times $k\leq r$, for all $r\geq2$,

  \item[(d)] Extremal associativity morphisms  $(\alpha')^{-1}$ and $(\beta')^{-1}$  coincide on all $h+1$-tuples of the form
 $(V, V, \dots, V, V^s)$ with $V$ repeated a number of times $h\leq s$, for all $s\geq2$.

\end{itemize}

  Then $\alpha=\beta$ and $c=d$ everywhere.
 
 \end{thm}
  
 \begin{proof}
 
 We modify the proof of Theorem \ref{claim0}.
 
 First of all, we either insert at the end of a paragraph x) or replace the whole paragraph x), with the paragraph x') as follows
 from a) to h).

After a) we insert

a')  Similarly, for the braiding morphisms: to show that
$c(V^r, V^s)=d(V^r, V^s)$ on a given pair of tensor powers of $V$, it suffices to assume that $\min\{r, s \}\geq 1$, as otherwise
the braiding morphisms equal the identity map, by definition.
\medskip

After b) we insert

 b') Similarly, on any pair $(r, s)$ of positive integers, we defined the  
 integer-valued function $$h(m, r):=m+r.$$ Thus $h(m, r)\geq 2$ and $h(m, r)=2$ if and only if $(m, r)=(1, 1)$;
   $h(m, r)=3$ if and only if $(m, r)=(1, 2)$ or $(m, r)=(2, 1)$.
  In these cases the two braiding morphisms $c$ and $d$ coincide on the corresponding pairs
    $(V, V)$, $(V, V^2)$  and $(V^2, V)$ by assumption as well.\medskip  
 We replace d) with 
 
 d') We proceed by induction. First inductive assumption on $f$ and $h$: 
\medskip

Let $N$ be an integer $\geq3$. Assume  that $\alpha_{V^s,V^t,V^u} = \beta_{V^s,V^t,V^u}$ for all triples $(s, t, u)$ such that $f(s, t, u) \leq N$. \medskip

Assume also that $c(V^m, V^r)=d(V^m, V^r)$ for all pairs $(m, r)$ such that $h(m, r)\leq N$.
\medskip

We replace e) with 

e')  First inductive step on $f$ and $h$: 

Let    
 $(V^s,V^t,V^u)$ be  a new triple with $$f(s,t,u) = N+1.$$
   Then at least one among $s$, $t$, $u$ takes a value $>1$. 
   \medskip
   
      Let $(V^m, V^r)$ be a new pair with $h(m, r)=N+1$. \medskip
      
      We repeat 
      
      f)  Second inductive assumption on $g$ 
with the value of both $f$ 
    fixed to $N+1$:
      
      The function $g$ takes minimum value $=2$, and this happens if and only if  two coordinates equal $1$ and the third is an arbitrary tensor power of $V$. On such triples, $\alpha$ and $\beta$ coincide by assumption.
Let $M$ be an integer $\geq2$.  
 It consists of the following 
 assumption on the associativity morphisms.

 Let us consider
  triples $(s, t, u)$ on which    $f$  takes the value $N+1$.
     Assume   that for each such $(s, t, u)$, we   have
    $$\alpha_{V^s,V^t,V^u} = \beta_{V^s,V^t,V^u}, \quad \text{for }
    g(s, t, u)\leq M,$$\medskip

 This assumption holds for the minimal value $M=2$ as before remarked.\medskip
 
 We repeat  
 
 g)  Second inductive step on $g$  
 (at $f=N+1$):
Let $(s, t, u)$ be a  new  triple for which $f(s, t, u)=N+1$
 and $g(s, t, u)=M+1$.
Then at least two coordinates among $(s, t, u)$ take a value $>1$.

\medskip

We replace h) with 

h') We first 
show that $\alpha=\beta$ in the second inductive step for $g=M+1$ (at $f=N+1$).
This is the main part of the proof. 

To complete the proof one also needs to show that equality of between the braiding morphisms $c$ and $d$
holds in  the first inductive step on the increased value of  $h=N+1$. This equality   will be deduced after the completion of     the proof of equality of associativity morphisms in the second inductive step on $g$ stated in the previous paragraph. See   the following paragraph r) for the braiding morphisms, where we shall use
  methods similar to those of the proof of Prop. 5.5 in  \cite{CGP}.
\medskip

To prove the statement in the first paragraph in h'), we first  notice that   paragraphs i), l), m), n) in the proof of Theorem 
\ref{claim0} hold also in the present case. 

Specifically,
  i) and l) only deal with associativity morphisms and pentagon equation.
  
  Moreover m)   only depends on naturality  of the $T$-morphisms on its three variables, and we apply this naturality to
the braided symmetry $c$ which we know to verify  the duality property for $V$. 

Equations in paragraph n)   now should be replaced by equations
  on the two hexagonal diagrams for the two structures $(\alpha, c)$ and $(\beta, d)$
  However, all
the braiding morphisms $d$  in the hexagonal equations used in n) for $(\beta, d)$     have the generating representation $V$ in one of the two variables.
By assumption, all these braiding morphisms in $d$ coincide with the corresponding braiding morphisms in $c$.
 Thus the equations in paragraph n) hold true also in the present   case.
\medskip

On the other hand, to obtain the equations in paragraph p) we have used   the hexagonal diagram (\ref{braided_symmetry1})  for 
$(\alpha, c)$ and $(\beta, d)$   with components on the pairs
$(V^u, V)$ on the horizontal right bottom arrow,  and also
on the pair
$(V^u, V^\ell)$ in the left vertical arrow
$(V^u, V^\ell\otimes V)$ in the  horizontal right upper arrow.
Similarly, for the use of  the hexagonal diagram (\ref{braided_symmetry2})
we need the components of $c$ and $d$ on the pairs $(V, V^u)$, $(V^\ell, V^u)$, $(V^\ell\otimes V, V^u)$.

By assumption we know that   $c$ and $d$ take the same value in the cases where $V$ appears as a component.
By the inductive assumption on the function $h$, we know that $c$ and $d$ coincide
on pairs $(V^u, V^\ell)$ and $(V^\ell, V^u)$, as $\ell+u=N$.

Thus to use paragraph p) in the present case we need to 
anticipate the proof that $c$ and $d$ coincide also on
$(V^u, V^\ell\otimes V)$ and $(V^\ell\otimes V, V^u)$.

We anticipate then the following paragraph
o), which shows  that $c$ and $d$ take the same
value also on such pairs.

o) We only show that
\begin{equation}
\label{equality_braiding_morphisms}
c(V^u, V^\ell\otimes V)=d(V^u, V^\ell\otimes V),
\end{equation} the other  case follow in a similar way. 

Recall that $\ell<u$,  and  that we are in   the second inductive step on $g$, so 
$f=\ell+u+1=N+1$ and
$g=\ell+1\leq M+1$.  

In the following passages we shall use   the induction hyothesis on $f$: 
$\alpha$ and $\beta$ coincide
on triples for which $f\leq N$,  

We shall also use the induction hypothesis
on  $c$ and $d$: they coincide on pairs for which $h\leq N$.

Moreover, we shall   use the assumptions (a), (b), (c), (d) in the statement.

\medskip

1) We first apply a right parenthesization to the   factor $V^\ell\otimes V$:  
$$V^u\otimes(V^\ell\otimes V)\longrightarrow
V^u\otimes( V\otimes(\dots\otimes(V\otimes V)\dots)).$$
Note that   on
  all the involved triples, $f$ takes a value $\leq\ell+1= N+1-u\leq N$, so the corresponding associativity morphisms
  $\alpha$ and $\beta$ coincide.
  
2)   Then we apply a complete left parenthesization  
  $$V^u\otimes(V\otimes(\dots\otimes(V\otimes V)\dots)) \longrightarrow (\dots(V^u\otimes V)\otimes\dots)\otimes V.$$
  This passage is the same for $\alpha$ and $\beta$ by our assumption (c)  in the statement of   Theorem \ref{claim1}.
  
3)   Then we apply a braiding morphism to the tensor product  $V^u\otimes V$,  
  $$(\dots(V^u\otimes V)\otimes\dots)\otimes V\longrightarrow (\dots(V\otimes V^u)\otimes\dots)\otimes V.$$ This
  is independent of the choice of the braiding morphisms $(\dots(c\otimes 1)\otimes\dots)\otimes 1$ or $(\dots(d\otimes 1)\otimes\dots)\otimes 1$ by our assumption when one of the variables is $V$.

4) Then  we apply associativity morphisms that act identically on the last factor $V$. These do not depend 
on the choice of $\alpha$ and $\beta$ as  $f$ takes value $\leq N$.
The associativity that we choose lead to an object of the form
$$(\dots(V\otimes V^u)\otimes\dots)\otimes V\longrightarrow (V\otimes (V^u\otimes V^{\ell-1}))\otimes V.$$

5) Then we apply a braiding morphism to $V^u\otimes V^{\ell-1}$. Since $h \leq N$, this braiding morphism
is the same for $c$ or $d$.
Thus we obtain the object
$$((V\otimes (V^u\otimes V^{\ell-1}))\otimes V\longrightarrow (V\otimes (V^{\ell-1}\otimes V^u ))\otimes V.$$

6) Then we reassociate again on all the factors leaving the last $V$ fixed, and this associativity morphism is the same for $\alpha$ and $\beta$
again  as $f\leq N$. We get to
$$(V\otimes (V^{\ell-1}\otimes V^u ))\otimes V\longrightarrow ((V\otimes V^{\ell-1})\otimes V^u )\otimes V.$$

7) Then we use the long path of the pentagon equation (\ref{pentagon_equation}), which depends on the triples
$(V, V^{\ell-1}, V^u)$, $(V, (V^{\ell-1}\otimes V^u), V)$, $(V^{\ell-1}, V^u, V)$. On all these triples
$\alpha$ and $\beta$ coincide, either  by $f\leq N$ (left and right triples) or  by assumption (a) (middle triple).
We get to the object
$$((V\otimes V^{\ell-1})\otimes V^u )\otimes V\longrightarrow V\otimes (V^{\ell-1}\otimes(V^u\otimes V)).$$

8) Then we use the braiding morphisms on $V^u\otimes V$. These again coincide for $c$ and $d$ since $V$ appears as one of the variables.
We get to  the object 
$$V\otimes (V^{\ell-1}\otimes(V^u\otimes V))\longrightarrow V\otimes (V^{\ell-1}\otimes(V\otimes V^u)).$$

9) Finally we apply in order a parenthesization on the right on all coordinates in $V$ except for the first factor on the left by  $f\leq N$,  
$$ V\otimes (V^{\ell-1}\otimes(V\otimes V^u))\longrightarrow V\otimes(V\otimes (V\otimes\dots (V\otimes V^u))\dots),$$

10) a parenthesization on the left on all coordinates $V$, $V, \dots$, $V$, $V^u$ by our assumption (d) in the statement of Theorem \ref{claim1}
and we obtain
$$V\otimes(V\otimes (V\otimes\dots (V\otimes V^u))\dots)\longrightarrow (\dots(V\otimes V )\dots )\otimes V)\otimes V^u$$

11) and a suitable parenthesization on all coordinates except for the last $V^u$ by   $f\leq N$
leading to the original object on the left:
$$(\dots(V\otimes V )\dots )\otimes V)\otimes V^u\longrightarrow (V^\ell\otimes V)\otimes V^u.$$

Composing these coinciding braiding and associativity morphisms in the two cases,
  we get   the desired braiding morphisms $$V^u\otimes(V^\ell\otimes V)\to (V^\ell\otimes V)\otimes V^u$$
for both the braided symmetries $c$ and $d$, which must coincide, and this proves (\ref{equality_braiding_morphisms}). 

We may now
repeat paragraphs p) and q) in the proof of Theorem \ref{claim0} to see that $\alpha=\beta$ on triples for which
$f=N+1$ and   $g=M+1$ completing the proof that $\alpha=\beta$ for $f=N+1$   (first paragraph in h')).\medskip

r) At the completion of the   induction on the function g we have also shown that   $c(V^u, V^{\ell+1})=d(V^u, V^{\ell+1})$ 
$c(V^{\ell+1}, V^u)=d(V^{\ell+1}), V^u)$   whenever $\ell+1\leq u$ by naturality of $c$ and $d$, and therefore
  for $h=N+1$   (second paragraph in h')).

\end{proof}

  \section{Duality for the fundamental representation (Proof of Th. \ref{Finkelberg_HL}(\lowercase{d}))
  }\label{12}\bigskip

 In this section, we provide the final step to complete the proof of part (d) of Theorem \ref{Finkelberg_HL}. The core of this reduction relies on the duality property between the truncated tensor powers of the fundamental representation and the representation of the braid group induced by the R-matrix.

  The first instance of  this property    is Schur-Weyl duality for the general linear group.
This duality says that the algebra of intertwining operators, 
from a tensor power of the vector representation $V$ of the general linear group
 to itself (also called centralizer algebras for the tensor powers of the vector representation) is generated by permutation operators  \cite{Weyl}.

 The generalization of Schur-Weyl duality  to classical Lie groups and quantum groups
   is
  at the core  of the theory of both simple classical Lie groups and quantum groups.
 It is an enormous subject   
 developed by  many authors, extending
   from the classical  more familiar cases of type $A$ and vector representation to other Lie types beyond $A$, with respect to a specified generating representation $V$, and to the  quantum group case.

 \subsection{Conclusion of the Proof}
  The uniqueness Theorems  \ref{claim0} and   \ref{claim1}  both apply under a central  duality
 property of the fundamental representation $V$ with respect to the braid group representation in
the quantum group tensor category.

  To apply these theorems, we need the property that the centralizer algebras $(V^{n}, V^{n})$ with respect to
  the action of the Drinfeld-Jimbo quantum group $U_q({\mathfrak g})$ are generated by the natural representation of the
  braid group induced by the $R$-matrix. 
 
 \subsection{Roots of unity and truncated tensor powers}
 In our case $q$ is a root of unity as usually assumed corresponding to the
  unitary cases of the associated fusion category, specified in \cite{Wenzl}, \cite{Sawin}, see also Def.  19.1, 20.2 in \cite{CGP}.
  The tensor powers $V^{n}$ appearing in our setting are canonical choices of truncations of the full tensor powers
  of the quantum group $U_q({\mathfrak g})$ obtained by induction over $n$. 
  
  The construction of the iterated truncations $V^n$ of $U_q({\mathfrak g})$ at the root of unity
is explained in \cite{Wenzl} and is the starting point to   construct our weak Hopf algebras $A({\mathfrak g}, q)$ in Sect.  32 of \cite{CGP}.

\subsection{Cohomological reduction via braid group duality}
 The duality  property between the {\it truncated tensor powers}
$V^n$  of the fundamental representation  
and the representation of the braid group in the morphism algebra $(V^n, V^n)$ is required
by our cohomological
theorems    \ref{claim0} and  \ref{claim1}. These theorems assure   coincidence
between  braiding morphisms (associativity morphisms resp.) on all pairs (triples resp.) of truncated tensor powers $V^n$--and thus
complete the proof of part (d) theorem \ref{Finkelberg_HL}.

In the application of  Theorems  \ref{claim0} and  \ref{claim1}, we
  work on the side of the quantum group fusion category and by the following results, our assumptions on braid group duality for the fundamental
  representation of   our weak  Hopf
 algebras $A_W({\mathfrak g}, q)$ constructed in Sect. 32 in \cite{CGP} are verified for the defining structure  arising from Tannakian construction on the quantum group fusion category
 ${\mathcal C}({\mathfrak g}, q)$ and Wenzl functor.

\subsection{Inductive decomposition and nonzero quantum dimension}
For all the Lie types,   the truncated tensor powers $V^{n}$   are  inductively defined over $n$ by a canonical decomposition, which is into irreducible components for ${\mathfrak g}\neq E_8$,
 of a full tensor product of representations $V_\lambda\otimes V$ of  $U_q({\mathfrak g})$, with $V_\lambda$ an arbitrary irreducible
 of the semisimplified fusion category.
Then one selects  irreducible summands with nonzero quantum dimension. 
 The fusion of the fundamental representation of $E_8$ is described in \cite{Wenzl} as well, and is more delicate.
 
 \subsection{Generating the centralizer algebras for specific Lie types}
 Returning to the question of whether the braid  group representation  generates the centralizer algebras $(V^{n}, V^{n})$
 for the truncated tensor powers of $U_q({\mathfrak g})$ at the root of unity, this question can be studied in the same way as for the case of a corresponding summand  
with $q$ not a root of unity, as in the truncated summand the root of unity property
is not visible. This is explained among other things in Prop. 2.4 in \cite{Wenzl}, see also Theorem  31.2 in \cite{CGP} for the same  result with a slightly more detailed proof, where an important role is played by the continuous curve
described  in that proposition, that passes through generic values of $q$. 
See also   Section 3.5 in \cite{Wenzl} for the special tensor products $V_\lambda\otimes V$ of representations
of $U_q({\mathfrak g})$. 
Thus it suffices to have the duality property between the full algebras $(V^{\otimes n}, V^{\otimes n})$ of  $U_q({\mathfrak g})$ when $q$
is not a root of unity.

We recall that the fundamental representation $V$ are listed in Sect. 3.5 in \cite{Wenzl} for all the Lie types,
and in particular are given by

\begin{itemize}
    \item \textbf{Type $A$ and $C$ (vector 
    representation):} The Lie types $A$ and $C$ for the vector representation are well known. If $q$ is not a root of unity, the  centralizer algebras 
$(V^{\otimes n}, V^{\otimes n})$ in the type $A$ case are described
by Hecke algebras, a quotient of the complex braid group algebra by well known generators and relations. The first result in this 
direction is due to   Jimbo \cite{Jimbo}.

When $V$ is the vector representation, for the  Lie type  $C$ (and also $B$, $D$ for this representation), 
Murakami and Birman and Wenzl introduced and studied
what are known as the BMW algebras  \cite{Murakami}, \cite{Birman_Wenzl}.
See also
\cite{GHJ}. 
Duality for the vector representation is of quantum groups was studied by Kirillov and Reshetikhin \cite{Kirillov_Reshetikhin}.
A summary of these results and references to the original papers for the vector representation for all the classical Lie types may be found in Theorem 10.2.5
of \cite{Chari_Pressley}.

    \item \textbf{Types $B$,  and $D$ (spinor representation):} In particular $V$ is the direct sum of the two spinor representations
  for type $D_n$ with $n$ even. The centralizer algebra of the fundamental representation for types $B$ and $D$ involve difficult
work by Wenzl and his coauthors. Central theorems are Theorem 5.2 in \cite{Wenzl_dualities} for the even case, and the previous
Theorem 3.3 in \cite{Wenzl_dualities0} for the odd case. A special case has been studied by Rowell and Wenzl in 
\cite{Rowell_Wenzl} in their study of a conjecture by Rowell and his coauthors regarding the characterization of finiteness of the image of the braid group in the morphism algebras
of tensor powers of a simple object of a braided fusion category by integrality of the squared Frobenius-Perron dimension of the object.

    \item \textbf{Type $G_2$ (the $7$-dimensional representation):} Type $G_2$ was shown by Lehrer and Zhang \cite{Lehrer_Zhang}, Morrison  \cite{Morrison}, and Martirosyan and Wenzl
 \cite{Martirosyan_Wenzl}.
 
 \end{itemize}

\section{Concluding remarks and outlook}\label{13}

This work provides a unified solution to two long-standing problems: the construction of quantum gauge groups for conformal nets, as proposed by Doplicher and Roberts, and the direct proof of the Finkelberg equivalence between quantum groups and vertex operator algebras, as sought by Huang. By constructing the global quantum gauge group $A_{W}(\mathfrak{g}, q)$ and the unitary Drinfeld twist $T$, we have demonstrated that the modular tensor category of a WZW model is intrinsically tied to the analytic structure of its underlying weak Hopf algebra $A_W({\mathfrak g}, q)$, which is a quotient of $U_q({\mathfrak g})$ as an algebra. Alternatively, one can consider the Zhu algebra $A(V_{{\mathfrak g}_k})$ with a weak quasi-Hopf algebra structure with $3$-coboundary associator, which is a quotient of $U({\mathfrak g})$ as an algebra.

\subsection{Unifying high and low dimensional AQFT}
Our results reconcile the ``classical'' Doplicher--Roberts program---originally developed for high-dimensional Algebraic Quantum Field Theory (AQFT)---with the low-dimensional setting of conformal nets and VOAs. In high dimensions, the symmetry is always a compact group; in low dimensions, the braiding dictates that the symmetry must be a quantum object. 

The construction of $A_{W}(\mathfrak{g}, q)$ realizes the Mack--Schomerus program by providing a concrete Hopf algebraic framework for these symmetries. The fact that the identification can be reduced to the fundamental representation $V_{\rm fund}$ highlights a deep structural simplicity: the centralizer algebras of the WZW model are ``large enough'' to be generated by the braid group image. This property, verified for classical types and $G_2$, ensures that the categorical ``gauge'' is uniquely fixed by the braiding data alone. \textit{Furthermore, the isometric nature of the twist $T$ establishes the transport of the unitary structure, a result that remained elusive in previous formal algebraic approaches.}

\subsection{Future directions: exceptional types and generalizations}

While this paper completes the program for the classical Lie types and $G_2$, several directions remain for future exploration:

\begin{itemize}
    \item \textbf{Exceptional Lie types:} Extending the uniqueness proof to the types $F_4$, $E_6$, $E_7$, and $E_8$. As noted in the technical development, these cases possess centralizer algebras that require further study regarding the existence of specific forbidden levels to confirm that the image of the braid group is sufficient to uniquely determine completely the associativity.

    \item \textbf{Quantum homogeneous spaces:} The framework developed here for $A_{W}(\mathfrak{g}, q)$ suggests a path toward a systematic study of quantum homogeneous spaces in the VOA context. This could lead to a broader classification of ``quantum subgroups'' acting as symmetries for VOA extensions or orbifolds.
\end{itemize}

In conclusion, by providing an explicit, analytic construction of the tensor equivalence \textit{mediated by the Radon-Nikodym property of the twist}, this work establishes the field-theoretic reconstruction of quantum symmetries, bypassing the \textit{rigidity obstructions, character-based verifications, and KZ equations}.
\appendix

\section{Primary fields, KZ equations, quantum groups and associated braided categories}\label{2}

\subsection{Primary field in conformal field theory}
Primary fields, introduced under a different name by Mack and Salam in the context of higher-dimensional conformal field theories  \cite{Mack-Salam}, were investigated in 2D by Belavin, Polyakov, and Zamolodchikov. They recognized that the local symmetry algebra is given by the infinite-dimensional Virasoro algebra  \cite{BPZ}, building upon previous work by Feigin and Fuchs \cite{Feigin_Fuks}  and Kac \cite{Kac_Virasoro} on its representation theory.

 Witten introduced an interacting $2D$ conformal field theory,   now known as the Wess-Zumino-Witten (WZW) model \cite{Witten_WZW}. He showed that this model possesses a symmetry satisfying the relations of an affine Kac-Moody algebra. In the same year, Knizhnik and Zamolodchikov showed that the correlation functions of a WZW model  satisfy a system of partial
 differential equations \cite{KZ}.  
This Kac-Moody framework was related to the Virasoro algebras by Goddard, Kent and Olive  \cite{GKO}, who introduced the coset construction to realize Virasoro minimal models directly from WZW models.

\subsection{Quantum groups and monodromy representations}
The mathematical formulation of the KZ equations was made rigorous by Tsuchiya and Kanie  \cite{Tsuchiya_Kanie}, who studied
the model for ${\hat{\mathfrak sl}}_2$ and
  showed
that the associated monodromy yields representations of the braid group that factor through   
the Hecke algebra at a root of unity. 

 The discovery of quantum groups first arose in physics literature in the work by Faddeev, introduced as noncommutative function algebras.
In the mid 80s breakthrough papers by Drinfeld and Jimbo discovered  them as quantization of universal enveloping
algebras, and developed their properties.
We refer to  \cite{Chari_Pressley} for an in-depth presentation and historical information.

Kohno proved that the monodromy representations of the KZ equations are equivalent to those derived from the universal $R$-matrix of a quantum group
 \cite{Kohno2}.
 Before recalling the subsequent deep developments of the theorem by Kohno,
  we need first review the emergence of tensor categories in   constructive conformal field theory.

\subsection{Modular tensor categories, their emergence  in   CFT and quantum groups}
Moore and Seiberg \cite{Moore-Seiberg1}, \cite{Moore-Seiberg2} derived
 certain equations from rational conformal field theory and   
discovered a connection between their equations and   the structure of a tensor category. 
The corresponding mathematical framework, specifically 
the notion  of ribbon braided tensor category, has since been  settled in mathematics,
we refer to \cite{EGNO}, \cite{ENO},  and references therein.

An independent derivation of strict braided tensor $C^*$-categories arose  from the axiomatic 
operator algebraic approach to QFT in the 1970s, see   \cite{DR1}, \cite{Mueger3}, \cite{Haag} and references therein.
The connections between tensor $C^*$-categories and compact quantum groups are studied  
\cite{CQGRC}.

  Turaev gave a definition of modular tensor category, based on previous work with Reshetikhin 
  on the construction of knot and 3-manifold invariants  see \cite{Turaev}. 
  Several authors   contributed to construct modular fusion categories form quantum groups
  at roots of unity, and   a more complete historical account on the quantum group side may be found  
   in \cite{Rowell2},  \cite{Sawin},
  see also   the introduction  of  \cite{CGP}. We denote by ${\mathcal C}({\mathfrak g}, q)$ 
  the resulting modular fusion category,
   with $q=e^{i\pi/d\ell}$, and $d\in\{1, 2, 3\}$ a suitable integer structure constant from Lie theory.
   Unitary tensor structures have been constructed by Wenzl and Xu \cite{Wenzl}, \cite{Xu_star}. 
    Wenzl's work  plays a major 
   role in our work, and we try to explain this in this paper.

\subsection{Quasi-Hopf algebras and the Drinfeld-Kohno theorem}
In the constructive approach to CFT, the Wess-Zumino-Witten model is of central interest.   Modular, or more generally ribbon rigid braided tensor category
  constructions involved many authors among physicists and mathematicians.
 A   main difficulty in constructing tensor categories from models of conformal field theories is the construction of the associativity morphisms compatible with naturally assigned braiding morphisms via the two hexagonal equations. Another main difficulty is to show rigidity of the category.

Non-strict   examples   motivated by   conformal field theory
arose from the work of Drinfeld on  quasi-Hopf algebras \cite{Drinfeld_quasi_hopf}.
 Let ${\mathfrak g}$ be a complex simple Lie algebra.   The mentioned problems   were first solved by Drinfeld  in 
\cite{Drinfeld_cocommutative}, \cite{Drinfeld_quasi_hopf} \cite{Drinfeld_galois} for
  the Drinfeld category with the introduction of the notion 
and construction of non-trivial examples of quasi-Hopf algebras over $U({\mathfrak g})[[h]]$, 
via the KZ differential equation. Drinfeld used this equation to construct the associativity
morphisms compatible with naturally assigned braiding morphisms. Furthermore he established rigidity of this category
by showing a remarkable ribbon, braided tensor
equivalence with the module category arising from the quantized universal enveloping algebra
$U_h({\mathfrak g})$, which is a Hopf algebra. To construct the ribbon, braided, tensor equivalence, he used Drinfeld coboundary matrix ${\overline R}$, which is
a modification of the $R$-matrix by a square root of the coboundary defined by the ribbon structure.

This is  
  his remarkable Drinfeld-Kohno theorem, expanding
the earlier theorem by Kohno  \cite{Kohno2}.
Drinfeld-Kohno theorem connects via a Drinfeld twist, the quasi-Hopf algebra $U({\mathfrak g})[[h]]$ with the
quantized universal enveloping algebra $U_h({\mathfrak g})$.
In particular,  
 the representations of the braid group arising from monodromy of the KZ equation are shown equivalent
to those arising from the action of the $R$-matrix of $U_h({\mathfrak g})$ with respect to the vector representation, 
for the classical Lie types.

\section{Affine Lie algebras, vertex operator algebras, the Kazhdan-Lusztig-Finkelberg equivalence,  the Verlinde formula, and the rigidity problem}\label{rigidity_problem}

  \subsection{Affine Lie algebras and the Finkelberg equivalence}

 The   work by Finkelberg is motivated by   work in physics by Moore and Seiberg \cite{Moore-Seiberg1}, \cite{Moore-Seiberg2} and  is based on the approach
  by Beilinson, Feigin,
   and Mazur  \cite{BFM} in algebraic geometry for the   braided tensor structure on $\tilde{\mathcal O}_{\ell}$.

    With $k$ a positive integer, defined as the physical level, and $\ell=k+\check{h}$, with  $\check{h}$ the dual Coxeter number of ${\mathfrak g}$, Finkelberg considered a certain semisimple subquotient $\tilde{\mathcal O}_{-\ell}$  of the representation category constructed by 
Kazhdan and Lusztig, paralleling the subquotient construction of the quantum group fusion category ${\mathcal C}({\mathfrak g}, q)$ in the setting of quantum groups at roots of unity. He defined an equivalence
$\tilde{\mathcal O}_{\ell}\to \tilde{\mathcal O}_{-\ell}$  
from  the category $\tilde{\mathcal O}_{\ell}$  of certain integral modules of the affine Lie algebra at      level $k$ introduced by Beilinson, Feigin, and Mazur
 to the subquotient rigid braided tensor category    $\tilde{\mathcal O}_{-\ell}$. 
 Combining with the equivalence by Kazhdan and Lusztig, Finkelberg derived
   an equivalence between   $\tilde{\mathcal O}_\ell \to {\mathcal C}({\mathfrak g}, q)$.
 See  Theorem 4.3 in \cite{Finkelberg}. 
     His work excluded   levels $E_6$, $k=1$; $E_7$, $k=1$; $E_8$, $k=1, 2$ because the same holds from Kazhdan-Lusztig work.

  By a private communication with Huang on the history of this problem, I learned that it   turned out that the construction  by Finkelberg described in \cite{Finkelberg} connecting
$\tilde{\mathcal O}_{\ell}$  with   Kazhdan and Lusztig structure inherited by $\tilde{\mathcal O}_{-\ell}$
could not suffice to derive rigidity of $\tilde{\mathcal O}_{\ell}$ from that of $\tilde{\mathcal O}_{-\ell}$.
This implied that the proof of the {\it tensor} equivalence $\tilde{\mathcal O}_{\ell}\to \tilde{\mathcal O}_{-\ell}$ was not complete in \cite{Finkelberg}, and was completed in \cite{Finkelberg_erratum}.

 \subsection{Rigidity in vertex operator algebra and affine Lie algebra   categories at positive integer levels}\bigskip
     As I learned from discussions with Huang, rigidity was not proved in  \cite{BFM}.
   Huang and Lepowsky     introduced and developed  a rigorous direct construction  of a modular (vertex) 
tensor category structure in a very general setting for  vertex operator algebras, including as a particular case the affine vertex operator algebras
at  positive integer levels $k$ that we denote by $V_{{\mathfrak g}_k}$. The construction has been completed by  Huang between 2005 and 2008   \cite{Huang(modularity)}, \cite{Huang2} following a general construction of tensor product theory
in the setting of vertex operator algebras by Huang and Lepowsky \cite{Huang_LepowskiI}, \cite{Huang_LepowskiII}, \cite{Huang_LepowskiIII}
and by Huang \cite{Huang1}, \cite{Huang_2}.
In particular  \cite{Huang1} gives the construction of associativity morphisms which is 
the first difficulty in the construction of the tensor structure. An historical account may be found in \cite{HL}. 
We also refer to the introductions of \cite{Huang_Lepowski_affine} and 
\cite{Huang2018}. These works give in particular the first mathematical proof
of modularity, 
that applies moreover to a large class of representation categories of vertex operator algebras.

In particular rigidity and modularity property in the setting of vertex operator algebras was first shown in \cite{Huang(modularity)}. 
To show modularity, and in particular rigidity, Huang first discovered the need of  the Verlinde formula.
The paper \cite{Huang_2} gives a proof of the Verlinde formula in the general case. See also \cite{Huang3} for the proof of
convergence and extension property
assumed in \cite{Huang1} and \cite{Huang4} for the proof of modular invariance conjecture
of Moore and Seiberg needed in the proof of the Verlinde formula.

  I have     been informed by Huang  that these doubts were motivated by two facts.
  On one hand,   his earlier direct communications  in the 90s had   clarified  
that rigidity was not proved in \cite{BFM}.  On the other, the mentioned property that his proof of rigidity in the setting of vertex operator
algebras involved the Verlinde formula indicated that more work was likely needed to complete the proof in \cite{Finkelberg}.
These    facts together with early claims of rigidity  in \cite{Finkelberg},   motivated 
    Huang to a   comparison of the results in \cite{Finkelberg}  with their results on rigidity and this led to  note a  gap in
 a proposition  of \cite{Finkelberg}
  with implications on the main result.

Huang's study had led him to identify the gap  in the fact that
  a suitable constant    needed to be shown to be nonzero in a proposition useful to obtain the equivalence between the two local systems considered by Finkelberg and Kazhdan and Lusztig, and this proposition  would imply  derivation of rigidity in his case from previous work.
On the other hand, Huang   proof of rigidity in the setting of vertex operator algebras
  had similarly  reduced to show nontriviality of a constant, and to prove this   he  had used
   the Verlinde formula as  above alluded to, which was a main result in their work. This had led Huang to
 inform Finkelberg about utility of  this formula  to show rigidity in his case.

As  later  reported by Finkelberg
in \cite{Finkelberg_erratum},  it  turned out that the completion    of the proof
of that proposition leading to rigidity of $\tilde{\mathcal O}_\ell$ and  the     equivalence $\tilde{\mathcal O}_\ell\to \tilde{\mathcal O}_{-\ell}$, 
involved the Verlinde formula as well (which   became another   main result in Finkelberg's paper).
This reflects also on Finkelberg composed equivalence   $\tilde{\mathcal O}_\ell \to
\tilde{\mathcal O}_{-\ell}\to {\mathcal C}({\mathfrak g}, q)$,
with Kazhdan-Lusztig
equivalence
 that moreover still passes through       the negative integer levels.
This implies in particular that the theorem does not include a few exceptional cases, including   ${\mathfrak g}=E_8$ and $k=2$, which correspond to cases
not treated in Kazhdan-Lusztig theory.

Following  the introductions in  \cite{Huang_Lepowski_affine}, \cite{Huang2018},    Finkelberg work  
 \cite{Finkelberg}, \cite{Finkelberg_erratum}  may be reinterpreted as  giving  the construction of a   braided tensor category structure on $\tilde{\mathcal O}_\ell$.
One may also conclude that these papers    give
a second proof of rigidity and modularity of this category.

An approach to modularity based on algebraic geometry may be found in Theorem 7.0.1 \cite{BK} but the proof of rigidity is not given and, as already mentioned it was known to experts that there was no proof based on their approach at the time. The author thanks Y.-Z. Huang for this information.
We point out the very recent   paper by Etingof and Penneys  where the authors show among other things that weak rigidity in the sense of \cite{BK}
is equivalent to rigidity for certain semisimple braided tensor categories of moderate growth with simple unit, see  \cite{EP} giving a third
proof    of rigidity   and give applications to a simplification of
Huang's proof of rigidity. They   also answer positively to a question posed in \cite{BK} that for every
finite semisimple category over an algebraically closed field the data of a modular functor
is essentially equivalent to the structure of a modular fusion category, completing the proof of rigidity in \cite{BK}.

On the other hand, Finkelberg equivalence  theorem is   stated as Theorem 7.0.2 in \cite{BK} without proof.
The braided tensor equivalence between $\tilde{\mathcal O}_\ell\to {\mathcal C}({\mathfrak g}, q)$  obtained from  the  combined works by Finkelberg, Kazhdan and Lusztig  
is considered indirect by several authors.

  \bigskip

\section*{Acknowledgments}
The author is grateful to Sergio Doplicher for proposing this problem and for the wonderful, inspiring scientific discussions during all the years she had the good fortune to share with him. 

She is grateful to Yi-Zhi Huang for several comments on the proof of modularity for fusion categories of vertex operator algebras and on the history of the Kazhdan--Lusztig--Finkelberg theorem. She is also grateful to Marco V. Giannone for pointing out the incompleteness of the proof of Lemma 7.10 in version 7 of \cite{CGP}, which led to the present Theorems \ref{claim0} and \ref{claim1}. 

She expresses her gratitude to Hans Wenzl for kind email correspondence regarding centralizer algebras for spinor representations, and for references on duality for tensor powers of the fundamental representation of Drinfeld--Jimbo quantum groups for the Lie types $B$, $D$, and $G_2$ that led to the present paper. She also thanks Makoto Yamashita for interesting discussions on \cite{GNY} and \cite{Wenzl_E_8} during his visit to the Department of Mathematics at Sapienza University of Rome, and for insightful comments on this work. Finally, she thanks Pavel Etingof for bringing \cite{EP} to her attention. 

The author is   grateful to an anonymous referee of an earlier stage of this project for constructive suggestions that   helped highlight the  advantages of this approach regarding uniform consideration of simply laced and non-simply laced Lie types.

The methods and results presented here are the outcome of a long-term project initiated in 2016. The author thanks Sapienza University of Rome for financial support through the \textit{Finanziamento di Ateneo per la Ricerca Scientifica} for the years 2018, 2019, 2020, 2021, and 2023.

   \end{document}